\documentclass[11pt,reqno]{article}
\usepackage[all]{xy}
\usepackage{amsmath}
\usepackage{amsthm}
\usepackage{amscd}
\usepackage{amssymb}

\input{epsf.tex}

\numberwithin{equation}{section}

\newcommand{\nb}[1]{#1\nobreakdash-}

\newcommand{\textmatrix}[4]{\bigl( \begin{smallmatrix} #1 & #2 \\ #3 & #4 
\end{smallmatrix} \bigr)}

\theoremstyle{definition}
\newtheorem*{Definition}{Definition}
\newtheorem*{Remark}{Remark}

\theoremstyle{plain}
\newtheorem{theorem}{Theorem}[section]
\newtheorem{proposition}[theorem]{Proposition}
\newtheorem{lemma}[theorem]{Lemma}
\newtheorem{corollary}[theorem]{Corollary}
\newtheorem{claim}[theorem]{Claim}

\newcounter{remarks}
{\paragraph*{Remarks}\smallskip
     \begin{list}{\arabic{remarks}. }{\usecounter{remarks}%
          \setlength{\leftmargin}{0in}%
          \setlength{\rightmargin}{0in}%
          \setlength{\labelsep}{0pt}%
          \setlength{\labelwidth}{0pt}%
          \setlength{\listparindent}{0pt}%
     }
}
{
\end{list}
}

\newcommand\ds\displaystyle
\newcommand\wt[1]{\widetilde{#1}}

\DeclareMathOperator{\Out}{Out}
\DeclareMathOperator{\Aff}{Aff}

\DeclareMathOperator{\Stretch}{Stretch}
\DeclareMathOperator{\rank}{rank}

\DeclareMathOperator{\SL}{SL}
\DeclareMathOperator{\GL}{GL}

\DeclareMathOperator{\PSL}{PSL}
\DeclareMathOperator{\Homeo}{Homeo}
\DeclareMathOperator{\QIMap}{\widehat{\QI}}

\DeclareMathOperator{\Isom}{Isom}

\DeclareMathOperator{\Aut}{Aut}
\DeclareMathOperator{\Stab}{Stab}

\DeclareMathOperator\nbhd{Nbhd}

\DeclareMathOperator\CAT{CAT}
\DeclareMathOperator\image{image}
\DeclareMathOperator\kernel{ker}

\DeclareMathOperator\Axis{Axis}
\DeclareMathOperator\Area{Area}
\DeclareMathOperator\VZ{VZ}

\newcommand\Isomplus{\Isom_{+}}
\newcommand\Affplus{\Aff_{+}}

\newcommand\Hull{\mathcal H}

\newcommand\R{{\mathbf R}}

\renewcommand\H{{\mathbf H}}
\newcommand\Hyp\H
\newcommand\hyp{\H}
\newcommand\Z{{\mathbf Z}}

\newcommand\solv{{\scshape solv}}
\newcommand\Solv\solv

\newcommand\inject{\hookrightarrow}
\newcommand\homeo{\approx}

\newcommand\infinity{\infty}
\newcommand\bndry{\partial}
\newcommand{\bdy}{\bndry}
\newcommand{\from}{\colon}
\def\composed{\circ}
\newcommand\suchthat{\bigm|}
\newcommand\inverse{{-1}}
\newcommand\inv{\inverse}

\newcommand\union{\cup}
\newcommand\Union{\bigcup}

\newcommand\absvalue[1]{{\left| #1 \right|}}
\newcommand\abs[1]{\absvalue{#1}}

\newcommand\generatedby[1]{{\langle #1 \rangle}}

\newcommand\Id{\text{Id}}

\newcommand\intersect{\cap}

\newcommand\U{\mathcal U}

\newcommand\Svarc{\v{S}varc}

\newcommand\subgroup{<}
\newcommand\semidirect{\rtimes}

\newcommand\Teichmuller{Teichm\"uller}
\newcommand\Poincare{Poincar\'e}
\newcommand\Mobius{M\"obius}
\newcommand\Mod{{\mathcal{M}}}
\newcommand\MCG{\Mod}

\DeclareMathOperator\QI{QI}
\newcommand\cross{\times}

\newcommand\Haus{{\mathcal H}}
\newcommand\G{{\mathcal G}}

\newcommand\coarsecap{\underset{c}{\cap}}
\renewcommand\O{{\mathcal O}}
\newcommand\M{M}
\renewcommand\P{{\mathbf P}}
\newcommand\cequiv{\underset{c}{\approx}}

\newcommand\fol{{\mathcal F}}
\newcommand\Fol{\fol}
\newcommand\F\fol

\newcommand\C{\mathcal C}
\newcommand\Teich{\mathcal T}
\newcommand\T\Teich

\newcommand\PMF{{\P\mathcal MF}}
\newcommand\MF{{\mathcal MF}}
\newcommand\PM{{\P\M}}

\DeclareMathOperator\Comm{Comm}
\DeclareMathOperator\QSym{QSym}

\DeclareMathOperator\Supp{supp}

\newcommand\ray[2]{\overrightarrow{[#1,#2)}}
\newcommand\geodesic[2]{\overleftrightarrow{(#1,#2)}}

\title{The geometry of surface-by-free groups}
\author{Benson Farb\thanks{Supported in part by NSF grant DMS 9704640
and by a Sloan Foundation Fellowship.}\ \ 
and Lee Mosher\thanks{Supported in part by NSF grant
DMS 9504946.}}

\begin{document}
\maketitle
\begin{abstract}
We show that every word hyperbolic, surface-by-(noncyclic) free group
$\Gamma$ is as rigid as possible: the quasi-isometry group of $\Gamma$
equals the abstract commensurator group $\Comm(\Gamma)$, which in turn
contains $\Gamma$ as a finite index subgroup.  As a corollary, two such
groups are quasi-isometric if and only if they are commensurable, and 
any finitely generated group quasi-isometric to $\Gamma$ must 
be weakly commensurable with $\Gamma$. We use quasi-isometries to compute
$\Comm(\Gamma)$ explicitly, an example of how quasi-isometries can actually
detect finite index information.  The proofs of these theorems involve
ideas from coarse topology,  Teichm\"{u}ller geometry,
pseudo-Anosov dynamics, and singular \Solv\ geometry.
\end{abstract}



\section{Introduction}

Let $\Sigma_g$ be a closed surface of genus $g\geq 2$, and let
$\MCG(\Sigma_g)=\pi_0(\Homeo(\Sigma_g) )$ denote the mapping class group
of $\Sigma_g$.  A \emph{Schottky subgroup} $H$ of $\MCG(\Sigma_g)$ is a 
free group of pseudo-Anosov mapping classes whose action on the
\Teichmuller\ space $\Teich(\Sigma_g)$ is ``weak convex cocompact''---the
group
$H$ has a limit set in Thurston's boundary of $\Teich(\Sigma_g)$, when
$\rank(H)\ge 2$ this limit set is a Cantor set, and the action of $H$ on the
``weak convex hull'' of the limit set is cocompact. Schottky subgroups of
$\MCG(\Sigma_g)$ exist in abundance: given any collection $\{\phi_1,\ldots
,\phi_r\}$ of pairwise independent pseudo-Anosov mapping classes of
$\Sigma_g$, for any sufficiently large positive integers $a_1,\ldots,a_r$
the elements $\{\phi_1^{a_1},\ldots ,\phi_r^{a_r}\}$ freely generate a
Schottky subgroup. See \S\ref{section:schottky} for a review of Schottky
groups, taken from \cite{Mosher:hypbyhyp} and \cite{FarbMosher:schottky}.

Given $g\geq 2$ and a Schottky subgroup $H<\MCG(\Sigma_g)\approx
\Out(\pi_1(\Sigma_g))$, one can construct a group $\Gamma_H$ given  by
the split extension 
$$1\to \pi_1(\Sigma_g)\to \Gamma_H\to H\to 1
$$
The groups $\Gamma_H$ are precisely the surface-by-free groups which 
are word hyperbolic \cite{FarbMosher:schottky}, and the construction of
\cite{Mosher:hypbyhyp} shows that they are abundant.

We are interested in studying quasi-isometries of the groups
$\Gamma_H$ for several reasons: the $\Gamma_H$ provide basic  examples of
rigidity theorems for word hyperbolic groups outside the  context of
negatively curved manifolds (see also \cite{Bourdon:BuildingQI},
\cite{KapovichKleiner:OneDBoundaries}); they are examples of groups which
can be viewed as phase spaces of dynamical systems arising from hyperbolic
endomorphisms of manifolds (see also \cite{FarbMosher:BSOne}, 
\cite{FarbMosher:BSTwo}, \cite{FarbMosher:ABC}); and they provide examples
of groups which are as rigid as possible in a very concrete sense (see
Theorem \ref{theorem:finiteindex} below). In particular, each
$\Gamma_H$ has finite index in its own abstract commensurator
$\Comm(\Gamma_H)$ as well as in its quasi-isometry group.  
The computation of $\Comm(\Gamma_H)$ given in Theorem 
\ref{theorem:comm} is explicit, one of the few instances outside lattices in
Lie groups where this has been done.

Finally, we propose the general problem of studying the asymptotic geometry
of extensions of surface groups $\Sigma$. These groups exhibit a beautiful 
and rich geometry which is encoded by a subgroup of $\MCG(\Sigma)$ acting
on the Teichm\"{u}ller space $\Teich(\Sigma)$. Some elements of the
geometry of such groups can be found in \cite{Mosher:hypext},
\cite{Mitra:trees} and \cite{FarbMosher:schottky}, but our needs will
require a somewhat involved account of this theory for extensions of
surface groups by Schottky groups, which we give in
\S\ref{section:geomtop}.

\subsection*{Statement of results}

Our first theorem gives a complete classification of the groups 
$\Gamma_H$ up to quasi-isometry. Although the usual quasi-isometry
invariants such as growth, ends, isoperimetric functions, etc., are the
same for all of these groups, they are only quasi-isometric under the
strictest algebraic conditions. 

Given a Schottky subgroup $H$ of $\MCG(\Sigma)$, we will show in
\S\ref{SectionSubcover} that there is a smallest \nb{2}orbifold
covered by $\Sigma$, denoted $\O_H$, such that $H$ descends via the
covering map $\Sigma \to\O_H$ to a Schottky subgroup of
$\MCG(\O_H)$, still denoted $H$. The orbifold $\O_H$ plays
an important role in describing the quasi-isometry class of the
surface-by-free group $\Gamma_H$.

\begin{theorem}[Classification Theorem]
\label{theorem:classification}
Given $g_1,g_2\geq 2$, let $H_1 \subgroup \MCG(\Sigma_{g_1})$ and $H_2
\subgroup \MCG(\Sigma_{g_2})$ by Schottky subgroups of rank $\ge 2$. The
following are equivalent:
\renewcommand{\labelenumi}{(\theenumi)}
\begin{enumerate}
\item
\label{ItemQI}
The surface-by-free groups $\Gamma_{H_1}$ and
$\Gamma_{H_2}$ are quasi-isometric.
\item
\label{ItemAbstractComm}
$\Gamma_{H_1}$ and $\Gamma_{H_2}$ are abstractly commensurable, meaning
that they have finite-index subgroups which are
isomorphic. 
\item
\label{ItemConcreteComm}
There is an isomorphism $\O_{H_1} \approx \O_{H_2}$ such that in the group
$\MCG(\O_{H_1}) = \MCG(\O_{H_2})$ the Schottky subgroups $H_1$ and $H_2$ are
commensurable, meaning that $H_1 \intersect H_2$ has finite index in each
of $H_1$ and $H_2$.
\item
\label{ItemSameLimitSet}
There is an isomorphism $\O_{H_1} \approx \O_{H_2}$ such that in the
group $\MCG(\O_{H_1}) = \MCG(\O_{H_2})$ the Schottky subgroups $H_1$ and
$H_2$ have the same limit set in the Thurston boundary of the
\Teichmuller\ space $\Teich(\O_{H_1}) = \Teich(\O_{H_2})$.
\end{enumerate}
\end{theorem}

The next theorem shows that each of the groups $\Gamma_H$ is determined
among all finitely-generated groups by its asymptotic geometry, up to
finite data.

\begin{theorem}[Quasi-isometric rigidity]
\label{theorem:rigidity} 
Let $\Gamma_H$ be a surface-by-free group with $H \subgroup
\MCG(\Sigma_{g})$ a Schottky group of rank $\ge 2$. If $G$ is any
finitely-generated group which is quasi-isometric to $\Gamma_H$, then there
is a finite normal subgroup $F \subgroup G$ such that $G/F$ is abstractly
commensurable to a surface-by-free group. Combining with Theorem
\ref{theorem:classification} it follows that $G/F$ is abstractly
commensurable to $\Gamma_H$. 
\end{theorem}

\paragraph{Remark.} We emphasize that it is essential for our methods that
$\rank(H) \ge 2$ in the statements of Theorems
\ref{theorem:classification} and \ref{theorem:rigidity}. Indeed, when $H$
is a rank~1 Schottky subgroup, in other words an infinite cyclic subgroup
generated by a pseudo-Anosov mapping class, then $\Gamma_H$ is the
fundamental group of a closed hyperbolic \nb{3}manifold that fibers over
the circle
\cite{Otal:fibered}, and hence all the groups $\Gamma_H$ with $H$ Schottky
of rank~1 are quasi-isometric to each other and to $\hyp^3$. Moreover, the
restatement of Theorem \ref{theorem:rigidity} (minus the last sentence)
for $H$ of rank~1 is equivalent to Thurston's virtual surface bundle
conjecture for closed  hyperbolic $3$-manifolds. \label{Footnote}
\footnote{This equivalence, however, hides the following fact: Theorem
\ref{theorem:classification} provides a commensuration from $\Gamma_{H_1}$
to $\Gamma_{H_2}$ coarsely taking $\pi_1(\Sigma_{g_2})$ to
$\pi_1(\Sigma_{g_2})$, a situation which is often impossible for rank~$1$;
see \S\ref{SectionHorResClasses}. Hence the proofs of Theorems
\ref{theorem:classification} and \ref{theorem:rigidity} cannot 
shed light on the virtual surface bundle conjecture.}

\bigskip

Theorem \ref{theorem:rigidity} and (most of) Theorem
\ref{theorem:classification} follow from our main result, Theorem
\ref{theorem:main} below.

\paragraph{Commensurations and quasi-isometries.}
Recall that a \emph{commensuration} of $\Gamma$ is an isomorphism between
finite index subgroups of $\Gamma$. Composition of two commensurations is
defined on a further finite index subgroup. Two commensurations are
equivalent if they agree on a common finite index subgroup.
Composition of equivalence classes gives a well-defined group operation,
and we thereby obtain the \emph{abstract commensurator group}
$\Comm(\Gamma)$ of the group $\Gamma$.  

The \emph{quasi-isometry group} $\QI(\Gamma)$ is the group of coarse
equivalence classes of self quasi-isometries of $\Gamma$ (endowed with
any word metric), where two quasi-isometries are coarsely equivalent if
they have finite distance in the sup norm.

In general it is difficult to compute the abstract commensurator of a
group. The computation for irreducible lattices in semisimple groups $G
\neq \PSL(2,\R)$ is the content of Mostow Rigidity together with theorems
of Borel and Margulis (see e.g.\ \cite{Zimmer:Book}). Our main result,
Theorem~\ref{TheoremCommQI}, says that the groups $\Comm(\Gamma_H)$ and
$\QI(\Gamma_H)$ are isomorphic, and gives an explicit computation of
these groups. This is the first time we know of where quasi-isometries
are used to compute the abstract commensurator group. While our
expression for $\Comm(\Gamma_H)$ is purely algebraic, we do not know how
to do the computation algebraically.

Recall that given groups $K \subgroup Q$, the
\emph{(relative) commensurator} of $K$ in $Q$, denoted $\Comm_Q(K)$, is
the subgroup of all $q\in Q$ such that conjugation by $q$ takes some
finite index subgroup of $K$ to another, or equivalently $K \intersect
qKq^\inv$ has finite index in both $K$ and $qKq^\inv$; this subgroup is
also known as the \emph{virtual normalizer} of $K$ in $Q$. The
group $\C=\Comm_{\MCG(\O_H)}(H)$, the relative commensurator of $H$ in
$\MCG(\O_H)$, plays a key role in computing the abstract commensurator
of~$\Gamma_H$:

\begin{theorem}[Computation of \protect{$\QI(\Gamma_H)$} and 
\protect{$\Comm(\Gamma_H)$}]
\label{TheoremCommQI}
\label{theorem:comm}
\label{theorem:main}
\label{theorem:finiteindex} 
Given $g\geq 2$ and a Schottky group $H \subgroup \MCG(\Sigma_g)$ of
rank $\ge 2$, the natural homomorphism 
$$\Comm(\Gamma_H) \to \QI(\Gamma_H) 
$$ is an isomorphism.  Furthermore, these groups are
isomorphic to the group $\Gamma_\C$ given 
explicitly by the short exact sequence 
$$1 \to \pi_1(\O_H) \to \Gamma_\C \to \Comm_{\Mod(\O_H)}(H) \to 1 
$$ 
Moreover, $\C=\Comm_{\Mod(\O_H)}(H)$ contains $H$ as a finite index 
subgroup. In particular, $\Gamma_H$ has finite index in
$\Comm(\Gamma_H)\approx\QI(\Gamma_H)$.
\end{theorem}

This theorem is proved in Sections~\ref{SectionComputingQIGroup}
and~\ref{section:comm}.

The fact that $\Gamma_H$ has finite index in $\QI(\Gamma_H)$ and also in 
$\Comm(\Gamma_H)$ shows that the groups $\Gamma_H$ are extremely rigid.
Among irreducible lattices in semisimple Lie groups, this phenomenon
holds for the noncocompact, nonarithmetic lattices  and for no other
lattice; see \cite{Margulis:DiscreteSubgroups} for the commensurator
statement, and \cite{Schwartz:RankOne} and e.g.\ \cite{Farb:Lattices} for
the quasi-isometry statement.

\subsection*{Outline of the proof of Theorem \ref{theorem:finiteindex}}

After some preliminary material, we begin in \S\ref{section:geomtop} by
constructing a quasi-isometric model space $X_H$ on which $\Gamma_H$
acts properly cocompactly by isometries.  A Schottky subgroup
$H<\MCG(\Sigma)$ gives an $H$-equivariant embedding of the Cayley graph
$T_H$ of $H$, a tree, in the Teichm\"{u}ller space
$\Teich(\Sigma)$. This embedding gives a $\Sigma$-bundle over $T_H$,
each fiber carrying a hyperbolic structure representing the point of
$T_H
\subset \Teich(\Sigma)$ over which that fiber lies. The universal cover of
this
$\Sigma$-bundle is our model space $X_H$, an $\hyp^2$-bundle over
the tree~$T_H$.

An isometry of $\Gamma_H$ acts on $X_H$, permuting or ``respecting''
various patterns of geometric objects, typically foliations. Indeed, the
same patterns respected by isometries of $X_H$ are also respected by the
finite index supergroup $\Gamma_\C$ which is defined by the short exact
sequence given in Theorem~\ref{TheoremCommQI}; see
\S\ref{SectionComputingQIGroup} for the precise definition of
$\Gamma_\C$. 

The rest of the proof is devoted to showing that an arbitrary
quasi-isometry $f\from X_H\to X_H$ coarsely respects so many patterns
that it must be close to an element of $\Gamma_\C$. Using coarse
topology, we show in \S\ref{section:geomtop} that $f$ permutes the
collection  of ``hyperplanes'' $P_w$; these are the
$\hyp^2$-bundles over bi-infinite lines $w$ in $T_H$. The Schottky
property is used to relate each line $w$ to a \Teichmuller\ geodesic,
which in turn allows us to impose extra structure on the hyperplane
$P_w$: a ``pseudo-Anosov flow'' and a singular \solv\ structure. In
\S\ref{section:dynamics} we use this structure to prove that $f$ coarsely
respects several dynamically defined foliations associated to these flows,
such as the stable and unstable foliations.

In \S\ref{section:periodic}, we apply the above together with R.\ Schwartz's
geodesic pattern rigidity \cite{Schwartz:Symmetric} to show that the
quasi-isometry $f$ actually permutes the collection of {\em periodic
hyperplanes}, i.e.\ those $P_w$ with $w$ a bi-infinite periodic line in
the tree $T_H$. 

In \S\ref{section:endgame} we study how the quasi-isometry $f$ acts on the
limit set of $H$ in Thurston's boundary of $\Teich(\Sigma)$, the space of
projective measured foliations. This limit set is a Cantor set, and the
set of periodic hyperplanes gives a countable dense subset which is
preserved by the action of $f$. This information is then used to show
that $f$ is close to an element of $\Gamma_\C$. 

The final parts of the proofs of all of the main theorems are contained in
\S\ref{section:endgame}.

\paragraph{Surfaces versus orbifolds.}

Our main theorems are stated solely for closed, oriented surfaces
$\Sigma_g$ of genus $g \ge 2$, and Schottky subgroups $H$ of
$\MCG(\Sigma_g)$. But the conclusions of these theorems force us to
consider a wider realm: closed \nb{2}orbifolds $\O$ and virtual Schottky
subgroups of $\MCG(\O)$; see e.g.\ Theorem \ref{theorem:comm}
in which, as the proof will show, $\Comm_{\MCG(\O_H)}(H)$ 
is a virtual Schottky subgroup of
$\MCG(\O_H)$. Various results that we will need about (virtual) Schottky
subgroups of orbifold mapping class groups are formulated in
\cite{FarbMosher:schottky}. 

This raises the question of whether the quasi-isometry classes of
(orbifold)-by-(virtual Schottky) groups constitute a wider universe than
the quasi-isometry classes of (surface)-by-(Schottky) groups. The answer is
no: these universes are identical. The proof is given in Section
\ref{SectionSurfaceVsOrbifold}.

\paragraph{Acknowledgements} We are grateful to Kevin Whyte for crucial
help in proving Theorem~\ref{TheoremCommQIIsomorphism}.


\section{Preliminaries}
In this section we briefly review some facts about Teichm\"{u}ller
space and about quasi-isometries.  For details
we refer to \cite{Abikoff:realanalytic}, \cite{FLP},
\cite{ImayoshiTaniguchi}.

For most of the paper, while we concentrate on the proof of
Theorem~\ref{theorem:finiteindex}, we shall fix the genus $g$ and denote
$\Sigma=\Sigma_g$.

\subsection{\Teichmuller\ space and mapping class groups}

The \emph{\Teichmuller\ space} of $\Sigma$, denoted
$\Teich=\Teich(\Sigma)$, has two equivalent descriptions, related to each
other by Riemann's uniformization theorem: $\Teich$ is the space of
conformal structures on $\Sigma$ modulo isotopy; or it is the space of
hyperbolic structures on $\Sigma$ modulo isotopy. A topology and a real
analytic structure on $\Teich$ is specified by the geodesic length
embedding $\Teich\to[0,\infinity)^\C$ where $\C$ is the set of isotopy
classes of nontrivial simple closed curves on $\Sigma$, and a hyperbolic
structure on $\Sigma$ determines an element of $[0,\infinity)^\C$ by
taking the length of the unique closed geodesic in each isotopy class.
With respect to this topology, $\Teich$ is homeomorphic to Euclidean
space of dimension $6g-6$.

The \emph{mapping class group} of $\Sigma$, denoted $\MCG=\MCG(\Sigma)$,
is the group $\pi_0(\Homeo(\Sigma)=\Homeo(\Sigma)/\Homeo_0(\Sigma)$ where
$\Homeo(\Sigma)$ is the group of all self-homeomorphisms of $\Sigma$ and
$\Homeo_0(\Sigma)$ is the normal subgroup of self-homeomorphisms isotopic
to the identity.

The group $\MCG$ acts real analytically on $\Teich$. Also, $\MCG$ acts on
$\C$ and so on $[0,\infinity)^\C$, and the embedding
$\Teich\to[0,\infinity)^\C$ is $\MCG$-equivariant. The action of $\MCG$ on
$\Teich$ is properly discontinuous. The quotient orbifold
$\M=\M(\Sigma)=\Teich/\MCG$ is called the \emph{moduli space} of
$\Sigma$. A subset of $\Teich$ is said to be \emph{cobounded} if its image
under the quotient map $\Teich\to\MCG$ is bounded.

\subsection{Measured foliations}
A \emph{measured foliation} on $\Sigma$ is a foliation on the complement
of a finite set of singularities, together with a positive, transverse
Borel measure, such that each singularity is an \emph{$n$-pronged
singularity} for some $n \ge 3$, locally modelled on the horizontal
measured foliation of the quadratic differential $z^{n-2} \, dz$ in the
complex plane. A \emph{saddle connection} of a measured foliation is leaf
segment which connects two distinct singularities, and \emph{Whitehead
equivalence} is the equivalence relation on the set of measured
foliations generated by isotopy and the collapse of saddle connections.
The \emph{measured foliation space} of $\Sigma$, denoted
$\MF=\MF(\Sigma)$, is the space of Whitehead equivalence classes of
measured foliations on $\Sigma$. A topology on $\MF$ is specified by the
``transverse measure embedding'' $\MF \to [0,\infinity)^\C$, where a
measured foliation determines an element of $[0,\infinity)^\C$ by taking
the infimum of the transverse measures of representatives of $\C$.

Given a measured foliation $\Fol$ and $r \in (0,\infinity)$, multiplying
the transverse measure on $\Fol$ by $r$ gives a new measured foliation
denoted $r\Fol$. This gives a free action of $(0,\infinity)$ on $\MF$,
whose quotient space is the space of projective measured foliations on
$\Sigma$, denoted $\PMF=\PMF(\Sigma)$. 

The embedding $\Teich \inject [0,\infinity)^\C$, composed with the
projectivization map $[0,\infinity)^\C \to \P[0,\infinity)^\C$, produces
an embedding $\Teich \inject \P[0,\infinity)^\C$. The embedding
$\MF\inject [0,\infinity)^\C$ induces an embedding
$\PMF\inject\P[0,\infinity)^\C$. Thurston's Compactification Theorem
\cite{FLP} says that image of $\overline\Teich=\Teich\union\PMF$ in
$\P[0,\infinity)^\C$ is a closed ball of dimension $6g-6$, whose interior
is $\Teich$ and whose boundary sphere is $\PMF$.

\subsection{Geodesics in $\Teich$} The \Teichmuller\ metric and its
geodesics are usually described in terms of holomorphic quadratic
differentials on Riemann surfaces. Using results of Gardiner and Masur
\cite{GardinerMasur} and of Hubbard and Masur \cite{HubbardMasur:qd}, the
metric can be presented directly in terms of measured foliations.

Consider of pair of measured foliations $\Fol_x,\Fol_y$ which are
\emph{transverse}, meaning that they have the same singular set, at each
singularity $\Fol_x$ and $\Fol_y$ have the same number of prongs, they are
transverse in the usual sense away from the singularities, and near an 
$n$-pronged singularity they are locally modelled on the horizontal and
vertical measured foliations of the quadratic differential $z^{n-2} dz$
on the complex plane. Let $\abs{dy}, \abs{dx}$ denote the transverse
measures on $\Fol_x,\Fol_y$, respectively; the leaves of $\Fol_x$
should be visualized as horizontal lines with transverse measure
$\abs{dy}$, and the leaves of $\Fol_y$ as
vertical lines with transverse measure $\abs{dx}$. The formula $dx^2 +
dy^2$ defines a singular Euclidean metric on $\Sigma$ denoted
$\mu(\Fol_x,\Fol_y)$, with total area
$$\Area(\mu)=\int_\Sigma \abs{dx} \, \abs{dy} < \infinity
$$
Underlying the metric $\mu(\Fol_x,\Fol_y)$ is a conformal structure on
the complement of the singularities, but the singularities are removable
and so we obtain a conformal structure on $\Sigma$ and a point in $\Teich$
denoted $\sigma(\Fol_x,\Fol_y)$. This gives a well-defined map
from a certain subset of $\MF\cross\MF$ to $\Teich$. Namely, letting $\U
\subset \MF\cross\MF$ denote the set of pairs $(\xi,\eta)$ which are
represented by a transverse pair $(\Fol_x,\Fol_y)$ of measured
foliations, it follows that $\sigma(\xi,\eta)=\sigma(\Fol_x,\Fol_y)$ is
well-defined independent of the choice of the representative pair
$(\Fol_x,\Fol_y)$, defining a map $\U \to \Teich$ (see \cite{GardinerMasur},
Theorem 3.1). 

A transverse pair $\Fol_x,\Fol_y$ is \emph{normalized} if
$\Area(\mu(\Fol_x,\Fol_y))=1$. Let $\U_0 \subset \U$ be the subset
represented by normalized transverse pairs. For each $(\Fol_x,\Fol_y) \in
\U_0$, the map $t \mapsto \gamma(t)=\sigma(e^t\Fol_x,e^{-t}\Fol_y)$ is a
real analytic embedding of $\R$ in $\Teich$; the image of this embedding
depends only on the projective classes $\xi = \P\Fol_x$, $\eta = \P\Fol_y$
and is denoted $\geodesic{\xi}{\eta}=\geodesic{\P\Fol_x}{\P\Fol_y}$. 
\Teichmuller's Theorem \cite{Abikoff:realanalytic} says that any two points
$p \ne q \in \Teich$ are contained in a unique such line
$\geodesic{\P\Fol_x}{\P\Fol_y}$; moreover, if $s,t \in \R$ are such that
$p=\sigma(e^t\Fol_x,e^{-t}\Fol_y)$ and
$q=\sigma(e^s\Fol_x,e^{-s}\Fol_y)$, then the formula
$$d(p,q) = \abs{s-t}
$$
gives a well-defined metric on $\Teich$, known as the \emph{\Teichmuller\
metric}. Each line $\geodesic{\xi}{\eta}$ then becomes a bi-infinite
geodesic in $\Teich$. If we restrict the parameterization $\gamma(t) =
\sigma(e^t\Fol_x,e^{-t}\Fol_y)$ to the half-line $t \in [0,\infinity)$
then we obtain a geodesic ray in $\Teich$ and we call the point $\P\Fol_x
\in \PMF$ the \emph{ending foliation} of the ray; if $\sigma=\gamma(0)$
and $\xi=\P\Fol_x$ then this ray is denoted $\ray{\sigma}{\xi}$. For any
\Teichmuller\ line $\geodesic{\xi}{\eta}$ the two points $\xi,\eta$ are
called the ending foliations of the line. Any two points $\sigma,\tau \in
\T$ are the endpoints of a unique finite geodesic segment, denoted
$\overline{\sigma\tau}$.

With respect to the \Teichmuller\ metric, $\Teich$ is a complete metric
space on which $\MCG$ acts by isometries. Royden's Theorem
\cite{Royden:AutIsomTeich} says, when the surface $\Sigma$ is closed and
oriented, that the homomorphism $\MCG \to \Isom(\Teich)$ is an isomorphism,
except for a small kernel on certain small surfaces: on $\Sigma_2$ the
single nontrivial element of the kernel being the hyperelliptic involution
of $\Sigma_2$.

Given $\sigma \in \T$ and $\xi \in \PMF$ there is exactly one ray in $\T$
with endpoint $\sigma$ and with ending foliation $\xi$; we denote
this ray $\ray{\sigma}{\xi}$. This gives a one-to-one correspondence
between $\T \cross \PMF$ and geodesic rays in $\T$.  Given $\xi,\eta\in
\PMF$, there exists at most one geodesic in $\T$ with ending foliations
$\xi,\eta$, and it exists if and only if $(\xi,\eta)$ is in the image of
$\U_0 \subset \MF\cross\MF$ under the projectivization map
$\MF\cross\MF\to\PMF\cross\PMF$. If this geodesic exists we
denote it $\geodesic{\xi}{\eta}$. This gives a one-to-one correspondence
between a certain subset of $\PMF\cross\PMF$ and the set of geodesics in
$\T$.

\begin{Remark} It is not in general true that the end of the ray
$\ray{\sigma}{\xi}$ converges in $\overline\T$ to the point $\xi$;
however, it is at least true for cobounded
rays. 
\end{Remark}

Let $T\T$ denote the tangent space of $\T$. There is an embedding
$\T\cross\PMF \mapsto T\T$, taking $(\sigma,\xi)$ to the tangent vector
at $\sigma$ of the ray $\ray{\sigma}{\xi}$ in $T\T$, denoted
$D\ray{\sigma}{\xi}$. The image of this embedding will be denoted
$T^1\T$, called the \emph{unit tangent bundle} of $\T$, and $T^1\T$ is,
in fact, a topological sphere bundle over $\T$. Moreover, the map
$\U_0\to T^1\T$ taking $(\xi,\eta)$ to $D\ray{\sigma(\xi,\eta)}{\eta}$, is
a homeomorphism. See \cite{GardinerMasur} and \cite{HubbardMasur:qd} for
proofs.

\begin{Remark} The \Teichmuller\ metric is not a Riemannian metric,
although it is a Finsler metric. As such, the unit tangent sphere
$T^1_\sigma\T$ at each $\sigma \in \T$ is not a true ellipsoid in the
vector space $T_\sigma\T$, but instead a more general convex, centrally
symmetric sphere \cite{Royden:AutIsomTeich}.
\end{Remark}

The flow on $\U_0$ defined by $(\xi,\eta) \cdot t = (e^t\xi,e^{-t} \eta)$
pushes forward under the homeomorphism $\U_0\to T^1\T$ to a flow on
$T^1\T$, namely the \emph{geodesic flow} of $\T$. In other words, each
tangent vector $v \in T^1\T$ is tangent to a unique geodesic $\gamma$,
and the geodesic flow $v \cdot t$ is obtained by pushing $v$ forward a
distance $t$ along $\gamma$.

\subsection{Pseudo-Anosov homeomorphisms}
A homeomorphism $h \from \Sigma \to \Sigma$ is \emph{pseudo-Anosov} if
there exists a transverse pair of measured foliations $\Fol_x,\Fol_y$ and
a $\lambda>1$ such that $h(\Fol_x)=\lambda\Fol_x$ and
$h(\Fol_y)=\lambda^\inv\Fol_y$; the foliations $\Fol_x,\Fol_y$ are called
the \emph{stable and unstable} measured foliations of $h$, and $\lambda$
is the \emph{expansion factor} of $h$. A mapping class $H \in \MCG$ is
said to be pseudo-Anosov if and only if it has a representative $h
\from\Sigma\to\Sigma$ which is pseudo-Anosov. 

By construction, a mapping class $H\in \MCG =\Isom(\Teich)$ is
pseudo-Anosov if and only if there exists a geodesic $\gamma$ in $\Teich$
such that $H(\gamma)=\gamma$ and the action of $H$ on $\gamma$ is a
translation of nonzero length. In this case, the geodesic $\gamma$ is
unique and is called the \emph{axis} of $H$, denoted $\Axis(H)$. Moreover,
if $h$ is a pseudo-Anosov homeomorphism representing $H$, with stable and
unstable foliations $\Fol_x,\Fol_y$ and expansion factor $\lambda$, then
$\Axis(H)=\gamma(\Fol_x,\Fol_y)$ and the translation length of $H$ equals
$\log(\lambda)$. 

Note that by a theorem of Bers \cite{Bers:ThurstonTheorem}, a mapping
class $H\in\MCG$ is pseudo-Anosov if and only if the function $\sigma
\mapsto d(H,H\sigma)$ has a positive minumum in $\T$; moreover this
minimum is acheived precisely on $\Axis(H)$.

\subsection{Quasi-isometries}

Given $K\geq 1,C\geq 0$, a $(K,C)$ \emph{quasi-isometry} between metric
spaces is a map $f\from X\to Y$ such that: 
\renewcommand{\labelenumi}{(\theenumi)}
\begin{enumerate}
\item For all $x_1, x_2 \in X$ we have
$$\frac{1}{K} \, d_X(x_1,x_2) - C \le d_Y(f(x_1),f(x_2)) \le K
d_X(x_1,x_2) + C
$$
\item $d_Y(y,f(X))\leq C$ for each $y\in Y$.
\end{enumerate}
If $f$ satisfies (1) but not necessarily (2) then it is called a 
$(K,C)$ \emph{quasi-isometric embedding}.  A quasi-isometric embedding 
$f:\R\to X$ is a {\em quasigeodesic} in $X$.  

A \emph{coarse inverse} of a quasi-isometry $f \from X \to Y$ is a
quasi-isometry $g \from Y \to X$ such that, for some constant $C'>0$, we
have $d(g\circ f(x),x)<C'$ and $d(f\circ g(y),y)<C'$ for all $x \in X$
and $y \in Y$. Every $(K,C)$ quasi-isometry $f \from X\to Y$ has a
$K,C'$ coarse inverse $g \from Y \to X$, where $C'$ depends only on
$K,C$: for each $y\in Y$ define $g(y)$ to be any point $x \in X$ such
that $d(f(x),y) \le C$.

A fundamental fact observed by Efremovich, by
Milnor \cite{Milnor:curvature}, and by \Svarc, which we use repeatedly
without mention, states that if a group $G$ acts properly
discontinuously and cocompactly by isometries on a proper geodesic
metric space $X$, then $G$ is finitely generated, and $X$ is
quasi-isometric to $G$ equipped with the word metric. 

Given a metric space $X$ and $A \subset X$, we denote $\nbhd_r(A)=\{x\in
X \suchthat d(x,r) \le A\}$, and given $A,B \subset X$, we denote the 
\emph{Hausdorff distance} by 
$$d_\Haus(A,B) = \inf \{r \in [0,\infinity] 
\suchthat A \subset \nbhd_r(B) \text{ and }B \subset \nbhd_r(A)\}
$$

Given a metric space $X$, the self quasi-isometries of $X$ are
denoted $\QIMap(X)$. Define the \emph{coarse equivalence} relation on
$\QIMap(X)$ where $f,g \in \QIMap(X)$ are coarsely equivalent, denoted $f
\cequiv g$, if 
$$\sup_{x\in X} d(fx,gx) = C < \infinity
$$
We call $C$ the
\emph{coarseness constant}. Composition of
elements of $\QIMap(X)$ gives a well-defined binary operation on the set of
coarse equivalence classes of self quasi-isometries of $X$, defining a
group $\QI(X)$, the \emph{quasi-isometry group} of $X$. If $h \from X \to
Y$ is a quasi-isometry of metric spaces then $h$ induces an isomorphism
$\QI(X)\to \QI(Y)$. In particular, when $\Gamma$ is a finitely generated
group the identity map is a quasi-isometry with respect to the word metrics
of any two finite generating sets, and so the quasi-isometry group
$\QI(\Gamma)$ is independent of choice of word metric on $\Gamma$.

A \emph{quasi-action} of a group $G$ on a metric space $X$ is a map $G
\cross X \to X$, denoted $(g,x) \mapsto g \cdot x$, such that for some $K
\ge 1, C \ge 0$ we have:
\begin{itemize}
\item For each $g \in G$ the map $x \mapsto g \cdot x$ is a
$(K,C)$-quasi-isometry.
\item For each $g,h \in G$, $x \in X$ we have $d(gh \cdot x, g \cdot (h
\cdot x)) \le C$; in other words, $L_{gh} \cequiv L_g \composed L_h$ with
coarseness constant independent of $g,h$, where $L_{\gamma}$ means left
multiplication by $\gamma$.
\end{itemize}
The quasi-action is \emph{cobounded} if there exists a bounded subset $D$
having nonempty intersection with every orbit of the quasi-action. The
quasi-action is \emph{proper} if for each $R>0$ there exists an integer
$m \ge 0$ such that for any $x,y \in X$ the cardinality of the set $\{g
\in G \suchthat (g \cdot B(x,R)) \intersect B(y,R) \ne \emptyset\}$ is at
most $m$.

A fundamental principle of geometric group theory says that if a finitely
generated group $G$ is quasi-isometric to a metric space $X$, then the
left action of $G$ on itself, when conjugated by a quasi-isometry $G \to
X$, defines a cobounded, proper quasi-action of $G$ on $X$. To be
precise, if we have coarsely inverse quasi-isometries $h \from X \to G$,
$\bar h \from G \to X$, then the formula $g \cdot x = \bar h(gh(x))$
defines a cobounded, proper quasi-action of $G$ on $X$.


\section{Schottky groups on Teichm\"{u}ller space}
\label{section:schottky}

In this section we recall from \cite{FarbMosher:schottky} the motivation
for and definition of Schottky subgroups of mapping class groups; see
that paper for details and proofs.

Recall that a \emph{Schottky group} in $\Isom(\H^n)$ is a discrete, free
{\marginpar \scriptsize use convex cocompact terminology}
subgroup $F$ such that every orbit is quasiconvex in $\Isom(\H^n)$.
Equivalently, letting $\Lambda$ be the limit set of $F$ and
$\Hull\Lambda$ the convex hull of $\Lambda$, the action of $F$ on
$\Hull\Lambda$ is cocompact; it follows that $\Hull\Lambda$ is
quasi-isometric to $F$, and this quasi-isometry extends continuously to
an $F$-equivariant homeomorphism between $\Lambda$ and the Gromov
boundary of $F$. This equivalence follows from the fact that $\Hyp^n$
itself is a $\delta$-hyperbolic metric space, and in fact the same (or
closely analogous) equivalence holds for free subgroups of word
hyperbolic groups, thereby providing a theory of Schottky subgroups of
word hyperbolic groups.

In \cite{FarbMosher:schottky} we mimic this setup for free subgroups $F$
of $\Isom(\Teich)=\MCG$. The tricky part is that \Teichmuller\ space
$\Teich$ is not $\delta$-hyperbolic. Nevertheless, when the properties
above for Schottky subgroups of $\Isom(\hyp^n)$ are carefully translated
into the language of $\Isom(\Teich)$, the results of Minsky
\cite{Minsky:quasiprojections} provide enough negative curvature in
$\Teich$ to prove the equivalence of various notions of Schottkiness for
$F$. In addition, one of the main  theorems of \cite{FarbMosher:schottky}
is that the Schottky condition  characterizes precisely those free
subgroups $F$ for which 
$\pi_1(\Sigma)\semidirect F$ is word hyperbolic.

\begin{theorem}[Schottky groups: Definitions]
\label{theorem:schottkydef}
Let $F$ be a finite rank free subgroup of $\MCG(\Sigma)$. The
following are equivalent:
\begin{description}
\item[1. Orbit quasiconvexity]
\label{ItemQuasiconvex}
Each orbit $\O$ of the action of $F$ on $\Teich$ is quasiconvex in
$\Teich$, i.e.\ there is a constant $A$ such that for any $x,y \in \O$, the
geodesic $\overline{xy}$ is contained in the $A$-neighborhood of $\O$.
\item[2. Weak convex cocompactness]
\label{ItemConvexcocompact}
There is a continuous, $F$-equivariant embedding of the Gromov
boundary $\bdy F$ into $\PMF$, with image denoted $\Lambda$, satisfying
the following:
\begin{enumerate}
\item For any $\xi \ne \eta \in \Lambda$ there is a geodesic
$\geodesic{\xi}{\eta}$ in $\Teich$; let
$$\Hull\Lambda =\union\Bigl\{\geodesic{\xi}{\eta} \suchthat
\xi\ne\eta\in\Lambda\Bigr\}
$$
be the \emph{weak convex hull} of $\Lambda$, and let
$\overline\Hull\Lambda =\Hull\Lambda\union\Lambda$.
\item
\label{ItemHypByHyp}
The $F$-equivariant homeomorphism $\Lambda \to \bdy F$ extends to
an $F$-equivariant map 
$$(\overline\Hull\Lambda,\Lambda,\Hull\Lambda) \to (F \union \bdy F, \bdy
F, F)
$$
which is continuous at each point of $\Lambda$ and which restricts to a
quasi-isometry $\Hull\Lambda\to F$, with respect to the \Teichmuller\
metric on $\Hull\Lambda$ and the word metric on $F$.
\end{enumerate}
\item[3. Word hyperbolic extension] The extension group
$\Gamma_F=\pi_1(\Sigma)\semidirect F$ is word hyperbolic.
\end{description}
\end{theorem}

\paragraph{Remark.} It follows from the weak convex cocompactness property
that for each geodesic $\geodesic{\xi}{\eta}$ in $\Hull\Lambda$, the image
of $\geodesic{\xi}{\eta}$ in $F$ is a quasigeodesic whose ends converge in
$F\union\bdy F$ to the images of $\xi,\eta$ respectively. It also follows
that each nontrivial element $f \in F$ is pseudo-Anosov, because $f$ has
an axis in the Cayley graph of $F$ and so $f$ has an axis in $\Teich$.

\bigskip

Theorem \ref{theorem:schottkydef} is proved in
\cite{FarbMosher:schottky}.  We call a subgroup $F$ satisfying any one of
the equivalent conditions of Theorem \ref{theorem:schottkydef} a
\emph{Schottky subgroup} of $\MCG(\Sigma)$, or a \emph{Schottky group} of
mapping classes.  These groups exist in abundance:

\begin{theorem}[Abundance of Schottky groups]
\label{theorem:existence}
Let $\phi_1,\ldots,\phi_n$ be a collection of $n$ independent 
pseudo-Anosov elements of $\MCG(\Sigma)$. Then for any sufficiently large 
natural numbers $a_1,\ldots,a_n$, the subgroup of $\MCG(\Sigma)$ generated by 
$\phi_1^{a_1},\ldots,\phi_n^{a_n}$ is a Schottky subgroup.
\end{theorem} 

\begin{proof} As noted in \cite{FarbMosher:schottky}, this follows from
Theorem \ref{theorem:schottkydef} together with the main result of
\cite{Mosher:hypbyhyp}.
\end{proof}


\section{The geometry and topology of $\Gamma_H$}
\label{section:geomtop}

\subsection{A geometric model for $\Gamma_H$}
 
Let $H$ be a Schottky subgroup of $\MCG=\Isom(\Teich)$.  We now
build a contractible, piecewise-Riemannian \nb{3}complex
$X_H$ on which $\Gamma_H$ acts freely, properly discontinuously, and
cocompactly by isometries, so that $\Gamma_H$ is quasi-isometric to $X$.

Choose a free generating set $h_1,\ldots,h_n \in\MCG$ for $H$.  
Let $\Delta$ be
a graph with $n$ edges
$\delta_1,\ldots,\delta_n$, each with one end at a common vertex $v_0$ of
valence $n$, and with opposite ends at valence one vertices
$v_1,\ldots,v_n$, respectively. Let $R$ be the rose with $n$ petals
obtained from $\Delta$ by identifying
$v_i$ with $v_0$ for each $i=1,\ldots,n$, and identify $\pi_1(R,v_0)$ with
$H$ so that the homotopy class of the loop $[\delta_i]$ is identified
with $h_i$.

On the product $\Sigma \cross \Delta$ make the following identifications:
for each $i$ choose a homeomorphism $\eta_i \from \Sigma\to \Sigma$
representing $h_i$ and identify $\Sigma \cross v_n$ with
$\Sigma \cross v_0$ by identifying $(x,v_0) \sim (\eta_i(x),v_i)$ for
each $x \in\Sigma$. Let $K_H$ be the quotient \nb{3}complex, and so we
obtain a locally trivial fiber bundle $K_H \to R$ with fiber $\Sigma$ and
with monodromy $H$. Up to a bundle isomorphism homotopic to the identity,
$K_H$ does not depend on the choices of representatives $\eta_i$. By Van
Kampen's Theorem we have an isomorphism $\pi_1(K_H)\approx \Gamma_H$. 

Let $X_H$ be the universal cover of $K_H$. Hence $\Gamma_H$ acts properly
discontinuously and cocompactly on $X_H$, with quotient $K_H$. We now
specify a metric on $X_H$ for which this action is isometric.

Let $T_H$ be the universal cover of $R$, an infinite, regular,
$2n$-valent tree, regarded as the Cayley graph for $H$. The universal
cover of $\Sigma$ is the \Poincare\ disc $D$. The bundle $K_H \to R$ with
fiber $\Sigma$ lifts to a locally trivial fiber bundle $\pi \from X_H \to
T_H$ with fiber $D$, and it follows that $X_H$ is homeomorphic to $D
\cross T_H$. In order to construct a $\Gamma_H$-equivariant metric on
$X_H$ the metrics on the fibers of $\pi$ must be appropriately
``twisted'' by considering the action of $H$ on $\Teich$.

We can embed the graph $\Delta$ in the tree $T_H$ as a fundamental
domain for the action of $H$, so that the restriction of the universal
covering map $T_H \to R$ agrees with the quotient map $\Delta \to
R$. Pick a base point $\sigma_0$ in the \Teichmuller\ space $\T$. Choose
a map $\rho \from \Delta \to \T$ taking $v_0$ to $\sigma_0$, taking
$v_i$ to $h_i(\sigma_0)$, and taking $\delta_i$ to a smooth path between
$\sigma_0$ and $h_i(\sigma_0)$, say the \Teichmuller\ geodesic. Extend
$H$-equivariantly to obtain a map $\rho \from T_H\to \T$. Henceforth we
shall fix the map $\rho$, and arcs in the tree $T_H$ will be
parameterized by arc length in $\Teich$ with respect to the map
$\rho$. We use $\tau$ as a variable taking values in $T_H$ and $d\tau$
as the arc length parameter in $T_H$.

Given an arc $\alpha \subset T_H$, we can impose a hyperbolic
structure on leaves of $\Sigma\cross \alpha$ in the form of a Riemannian
metric $g_\tau$ on each $\Sigma\cross \tau$, so that the conformal class
of $g_\tau$ represents the point $\tau \in \T$, and so that the
metrics $g_\tau$ vary smoothly with $\tau\in \alpha$. This can be obtained
for instance by choosing a pair-of-pants decomposition of $\Sigma$ and
using the associated Fenchel-Nielsen coordinates for $\T$ (see
\cite{FarbMosher:schottky} for a further discussion).
Then we may extend to obtain a piecewise smooth Riemannian metric on
$\Sigma\cross\alpha$ by the formula
$$ds^2 = g_\tau^2 + d\tau^2
$$
This metric is smooth on $\Sigma\cross\alpha$ except over the vertices of
$T_H$ in $\alpha$.

Applying this to each edge $\delta_i$ of $\Delta$ we obtain hyperbolic
structures on the fibers of $\Sigma \cross \Delta$, varying smoothly over
each edge $\delta_i$, so that for each $\tau \in\Delta$ the hyperbolic
surface $\Sigma\cross \tau$ represents $\rho(\tau) \in \T$. We therefore
obtain a piecewise smooth Riemannian metric on the \nb{3}complex
$\Sigma\cross\Delta$. Since $h_i(v_0)=v_i$ it follows that there is a
unique isometry $\eta_i\from\Sigma \cross v_0\to\Sigma\cross v_i$
representing the mapping class $h_i$, and we use these choices of $\eta_i$
to construct $K_H$ (here, for each $\tau\in\Delta$, we are implicitly
identifying $\Sigma\cross \tau$ with $\Sigma$ by projection to the first
factor). Thus we have defined a piecewise smooth Riemannian metric on
$K_H$ whose restriction to each fiber of the bundle $K_H\to R$ is a
hyperbolic metric on that fiber.

Lifting the metric on $K_H$ we obtain a piecewise Riemannian metric on
the universal cover $X_H$, whose restrictions to the fibers of the
bundle $\pi \from X_H \to T_H$ give a continuously varying family of
hyperbolic metrics $\overline g_\tau$ on the fibers, parameterized by
$\tau
\in T_H$. Each fiber is isometric to $\hyp^2$. For each arc $\alpha$ of
$T_H$ the sub-bundle $\pi^\inv(\alpha)$ has the form $D \cross \alpha$ and
the Riemannian metric has the form 
$$d \overline s^2 = \overline g_\tau^2 + d\tau^2
$$
where $\overline g_\tau$ is the lift to $D \cross \tau$ of the metric
$g_\tau$ on the appropriate fiber of $K_H$.

With respect to the piecewise Riemannian metric and the associated
geodesic metric on $X_H$, the group $\Gamma_H$ acts properly
discontinuously and cocompactly by isometries.

\subsection{What it means to coarsely respect a pattern}
\label{SectionPatterns}

A \emph{pattern} in a metric space is simply a collection of subsets. We
will loosely use the term \emph{foliation} to refer to a pattern forming
a partition of the space into disjoint subsets called \emph{leaves}, and
in that context the pattern itself will be called the \emph{leaf space}.

Let $X,Y$ be metric spaces, and let $\F,\G$ be patterns in $X,Y$
respectively. A quasi-isometry $\phi\from X \to Y$ is said to
\emph{coarsely respect} the patterns $\F,\G$ if there exists a number $A
\ge 0$ and a map $h \from \F\to\G$ such that for each element $L \in \F$
we have
$$d_\Haus(\phi(L),h(L)) \le A
$$ 
and if further a similar statement holds for a coarse inverse of
$\phi$. When distinct elements of $\F$ have infinite Hausdorff distance
in $X$, and similarly for $\G$, then $h$ is a bijection and it is uniquely
determined by the quasi-isometry $\phi$.

An isometry of $X_H$ respects many different patterns in $X_H$. The idea
of the proof of Theorem~\ref{theorem:finiteindex} is to show that an
arbitrary quasi-isometry $f$ of $\Gamma_H$ (hence of $X_H$) preserves
more and more structure, in the sense of coarsely respecting finer and
finer patterns, until so much structure is preserved that $f$ must
actually be a bounded distance from an element of the extension group
$\Gamma_\C$ defined in Theorem~\ref{TheoremCommQI}.

\subsection{Hyperplanes and the horizontal foliation}

The bundle $\pi \from X_H \to T_H$ associated to the group $\Gamma_H$
determines two important patterns in $X_H$: the \emph{horizontal
foliation} and the pattern of \emph{hyperplanes}.

The \emph{horizontal foliation} of $X_H$ is the pattern of fibers of
$\pi$, subsets $D_\tau=\pi^\inv(\tau)$, $\tau \in T_H$, each called a
\emph{horizontal leaf} of $X_H$. The leaf space of the horizontal
foliation is identified with $T_H$ via the bundle map $\pi$. Any two
horizontal leaves have finite Hausdorff distance in $X_H$, in fact from
the form of the metric on $X_H$ we see easily that $d_\Haus(D_\tau,D_u) =
d(\tau,u)$ for all $\tau,u \in T_H$. In other words, the bundle map $\pi$
induces an isometry between the horizontal foliation equipped with the
Hausdorff metric and the metric tree $T_H$. As noted earlier, the metric
on $X_H$ restricts to a metric $\overline g_\tau$ on each horizontal leaf
$D_\tau$ making $D_\tau$ isometric to $\hyp^2$.

A \emph{hyperplane} in $X_H$ is any set $P_{w}=\pi^\inv(w)$, where $w$ is
a bi-infinite geodesic in $T_H$. The metric on $X_H$ restricts to a 
piecewise Riemannian metric on $P_w$ of the form $\tilde g_\tau^2 +
d\tau^2$, with respect to the coordinates $P_w \homeo D \cross w$. Any
two distinct hyperplanes in $X_H$ have infinite Hausdorff distance. 

The pattern of hyperplanes in $X_H$ is the first pattern which must be
coarsely respected by a quasi-isometry.

\begin{proposition}[Hyperplanes respected]
\label{proposition:hyperplanes:respected}
Any quasi-isometry  $f\from X_H\to X_H$ coarsely respects the 
pattern of hyperplanes in $X_H$, with coarseness constant depending only on
the quasi-isometry constants of $f$.
\end{proposition}

In other words, $f$ must map each hyperplane in $X_H$ a uniform 
Hausdorff distance from some other hyperplane in $X_H$. Since distinct
hyperplanes have infinite Hausdorff distance in $X_H$ it follows that $f$
induces a bijection on the pattern of hyperplanes; when $f$ is understood
we denote this bijection $P_w\mapsto P_{w'}$.

Each hyperplane $P_w$ is uniformly properly embedded in $X_H$, that is,
there is a proper function $r \from [0,\infinity) \to [0,\infinity)$
independent of $w$ such that for all $x,y \in P_w$ we have:
$$d_{X_H}(x,y) \ge r\bigl( d_{P_w}(x,y)\bigr)
$$
By applying Proposition \ref{proposition:hyperplanes:respected} together
with a simple general principle (made explicit in Lemma 2.1 of
\cite{FarbMosher:ABC}) we have for any quasi-isometry 
$f:X_H\to X_H$ the following fact: for each hyperplane $P_w$, the map 
$f$ induces via restriction composed with
nearest-point projection, a quasi-isometry $\phi\from P_w\to P_{w'}$; the
quasi-isometry constants for $\phi$ depend only on those for
$f$.

The horizontal foliation on $X_H$ restricts to a horizontal foliation on
$P_w$. Using the coordinates $P_w \homeo D\times \R$ as
described above, the horizontal leaves in $P_w$ are of the form $D_t = D
\cross t$, $t \in \R$. The leaf space of this foliation is
$\R$ and we have $d_\Haus(D_s,D_t) = \abs{s-t}$ using Hausdorff distance
in $P_w$.

The fact that $X_H$ fibers over a tree $T_H$ with nontrivial branching 
restricts the behavior of $\phi$ as follows: 
 
\begin{proposition}[Horizontal foliation respected]
\label{proposition:horizontal:respected} For any quasi-isometry $f \from
X_H \to X_H$ and any hyperplane $P_w$, the induced quasi-isometry
$\phi\from P_w\to P_{w'}$ uniformly coarsely respects horizontal
foliations.
\end{proposition}

In other words, there exists a constant $A$ \emph{independent of $w$} such
that for any horizontal leaf $D_t \subset P_w$ there is a horizontal leaf
$D_{t'} \subset P_{w'}$ with $d_\Haus(\phi(D_t),D_{t'}) \le A$

This is precisely where the assumption is used that $H$ has rank greater
than one.

\begin{proof}[Proof of Propositions
\ref{proposition:hyperplanes:respected} and
\ref{proposition:horizontal:respected}]

Since $H$ has rank greater than one, its Cayley graph $T_H$ is a
\emph{bushy} tree, meaning that each point of $T_H$ is within some fixed
distance $\beta=1$ of some vertex $v$ such that $T-v$ has at least 3
unbounded components.  We can thus apply the following result, which is
Theorem 7.7 of \cite{FarbMosher:ABC}, to the metric fibration $\pi\from
X_H\to T_H$.

\begin{lemma}
\label{PropABCCoarseRespect}
Let $\pi \from X \to T$, $\pi' \from X' \to T'$ be metric fibrations
over bushy trees $T,T'$, such that the fibers of $\pi$ and $\pi'$
are contractible $n$-manifolds for some $n$. Let $f\from X \to X'$ be a
quasi-isometry. Then there exists a constant $A$, depending only on the
metric fibration data of $\pi,\pi'$, the quasi-isometry data for $f$, and
the constant $\beta$, such that:
\begin{itemize}
\item[(1)] For each hyperplane $P \subset X$ there exists a unique
hyperplane $Q \subset X'$ such that $d_\Haus(f(P),Q) \le A$.
\item[(2)] For each horizontal leaf $L \subset X$ there is a horizontal
leaf $L' \subset X'$ such that $d_\Haus(f(L),L')
\le A$.
\end{itemize}
\end{lemma}

Proposition \ref{proposition:hyperplanes:respected} is an immediate
consequence. To obtain Proposition \ref{proposition:horizontal:respected},
consider a horizontal leaf $D_t$ of $P_w$. By applying Lemma 
\ref{PropABCCoarseRespect} we obtain a horizontal leaf $D_s$ of $X_H$
uniformly Hausdorff close to $f(D_t)$. But $f(P_w)$ is uniformly Hausdorff
close to $P_{w'}$ and so $D_s$ is uniformly Hausdorff close to some
horizontal leaf $D_{t'} \subset P_{w'}$. Finally, the closest point
projection $f(P_w) \mapsto P_{w'}$ moves points a uniformly bounded
distance, and so $f(P_w)$ is a uniformly Hausdorff close to its closest
point projection in $P_{w'}$.
\end{proof}

\subsection{The singular \Solv\ metric on a hyperplane}
\label{section:singularSolv}

In this subsection we find a different metric on hyperplanes $P_w$ that
will allow us  to apply ideas from pseudo-Anosov
dynamics.

Associated to each geodesic $\gamma$ in $\T$ we recall the construction
of a ``singular \Solv\ metric'' on $D \cross \gamma$. Choose a
normalized transverse pair of measured foliations $\Fol_x, \Fol_y$ on
$\Sigma$, with transverse measures $\abs{dy}, \abs{dx}$ respectively, so
that $\gamma=\gamma(\Fol_x,\Fol_y)$; we can parameterize $\gamma$
by arc length:
$$\gamma(t) = \sigma(e^t \Fol_x, e^{-t} \Fol_y)
$$
For each $t$ we have a singular Euclidean metric
$\mu_t=\mu(e^t\Fol_x,e^{-t}\Fol_y)$, that is,
$$d\mu_t^2 = e^{-2t} dx^2 + e^{2t} dy^2
$$
Note that the transverse pair $(e^t\Fol_x,e^{-t}\Fol_y)$ is normalized
for all $t$, that is, $\Area(\mu_t)= \int_\Sigma e^t \abs{dy} e^{-t}
\abs{dx} = \int_\Sigma \abs{dx} \, \abs{dy} = 1$. Lifting to the
universal cover
$D=\wt\Sigma$ we obtain a transverse pair of measured foliations
$\wt\Fol_x,\wt\Fol_y$ with transverse measures still
denoted $\abs{dy},\abs{dx}$, and a one-parameter family of singular
Euclidean metrics $\tilde\mu_t$ on $D$, given by
$$d\tilde\mu_t^2 = e^{-2t} dx^2 + e^{2t} dy^2
$$
Now we can define a singular Riemannian metric on $D \cross \gamma$ by the
formula
\begin{align*}
ds^2 &= d\tilde\mu_t^2 + dt^2 \\
 &= e^{-2t} dx^2 + e^{2t} dy^2 + dt^2
\end{align*}
The singular locus of this metric is a family of vertical lines in $D
\cross \gamma$, whose intersection with the fiber $D \cross t$ is
precisely the singular set of $e^t\wt\Fol_x, e^{-t}\wt\Fol_y$. The
Riemannian metric extends across the singular locus to a complete
geodesic metric on
$D\cross \gamma$. The space $D\cross\gamma$ with this metric is called the
\emph{singular \solv\ space} associated to $\gamma$, denoted
$Q_\gamma$. The reason for the terminology is that away from the singular
locus the metric is locally modelled on \nb{3}dimensional
\solv-geometry, and at the singular locus (each component of which is a
line) the metric is modelled on some number of ``half \solv 's'' glued
together. 

Note that the singular \solv-metric on $Q_\gamma$ is uniquely
determined by the fibration $Q_\gamma\to\gamma$ together with the family
of $xy$-structures $e^t\wt\Fol_x, e^{-t}\wt\Fol_y$: in each leaf; the
measured foliation $e^t\wt\Fol_x$ provides the $e^{2t}dy$ term and the
measured foliation $e^{-t}\wt\Fol_y$ provides the $e^{-2t}dx$ term, and
the map to $\gamma$ with the arc length parameter $t$ provides the $dt^2$
term. 

In summary, we have a bijective correspondence
$$\{\text{geodesics $\gamma$ in $\T$}\} \leftrightarrow
\{\pi_1(\Sigma)\text{-equivariant singular \Solv\ metrics on $D\times
\R$}\}
$$

\paragraph{$xy$-structures.} For convenience we shall define an
\emph{$xy$-structure} on a surface
$S$ to be a transverse pair of measured foliations $\Fol_x,\Fol_y$. We will
work with $xy$-structures on $\Sigma$ as well as the lifted structures on
the universal cover~$D$. 

An \emph{affine automorphism} of an $xy$-structure
on $S$ is a homeomorphism of~$S$ which respects the (unordered) pair of
(unmeasured) foliations $\{\Fol_x,\Fol_y\}$, multiplying the transverse
measure on one by $\lambda>0$ and multiplying the transverse measure on the
other by $1/\lambda$. Let $\Aff(\Fol_x,\Fol_y)$ be the group of affine
automorphisms of the $xy$-structure $\Fol_x,\Fol_y$. There is a subgroup of
index $\le 2$ in $\Aff(\Fol_x,\Fol_y)$, denoted $\Affplus(\Fol_x,\Fol_y)$,
which respects the \emph{ordered} pair of (unmeasured) measured foliations
$(\Fol_x,\Fol_y)$, with transverse measures multiplied as described above.
There is a homomorphism 
$$\Stretch \from \Aff(\Fol_x,\Fol_y) \to \SL(2,\R)
$$
whose image lies in the subgroup of $\SL(2,\R)$ corresponding to matrices
which are either zero off the diagonal or zero on the diagonal, as follows.
Given $\phi\in \Aff(\Fol_x,\Fol_y)$: if $\phi(\Fol_x)=\lambda\Fol_x$ and
$\phi(\Fol_y)=\frac{1}{\lambda} \Fol_y$ then $\Stretch(\phi)$ is the matrix
$\textmatrix{\lambda}{0}{0}{1/\lambda}$; whereas if
$\phi(\Fol_x)=\lambda\Fol_y$ and
$\phi(\Fol_y)=\frac{1}{\lambda}\Fol_x$ then $\Stretch(\phi)$ is the matrix
$\textmatrix{0}{\lambda}{1/\lambda}{0}$. In the cases of interest where
$S=\Sigma$ or $D$, the image of the homomorphism, denoted
$\Stretch(\Fol_x,\Fol_y)$, is a discrete subgroup of $\SL(2,\R)$;
discreteness follows easily using the fact that the singular set of the pair
$\Fol_x,\Fol_y$ forms a nonempty, discrete net in $S$. As a consequence,
$\Stretch(\Fol_x,\Fol_y)$ is isomorphic to either $D_\infinity$, $\Z$,
$\Z/2$, or the trivial group. The kernel of the stretch homomorphism is the
group $\Isomplus(\Fol_x,\Fol_y)$ of automorphisms of the
\emph{ordered} pair of \emph{measured} foliations $\Fol_x,\Fol_y$.

We obtain a short exact sequence
$$1 \to \Isomplus(\Fol_x,\Fol_y) \to \Aff(\Fol_x,\Fol_y) \to
\Stretch(\Fol_x,\Fol_y) \to 1
$$
where $\Isomplus(\Fol_x,\Fol_y)$ is the group of automorphisms of the
\emph{ordered} pair of \emph{measured} foliations $\Fol_x,\Fol_y$.

Given a geodesic $\gamma=\gamma(\Fol_x,\Fol_y)$ in $\T$,
the universal cover $D=\wt\Sigma$ with the lifted $xy$-structure
$\wt\Fol_x,\wt\Fol_y$ may be identified with a certain horizontal fiber of
$Q_\gamma$. With respect to this identification, every isometry of the
singular \solv-manifold $Q_\gamma$ respects the horizontal foliation, and
vertical projection of an isometry onto $D$ defines an $xy$-affine
automorphism of $\wt\Fol_x,\wt\Fol_y$. Every $xy$-affine automorphism arises
in this manner, leading to an isomorphism between $\Isom(Q_\gamma)$ and
$\Aff(\wt\Fol_x,\wt\Fol_y)$. This leads in turn to an isomorphism of short
exact sequences as follows:
$$
\xymatrix{
1 \ar[r] & \Isomplus(\wt\Fol_x,\wt\Fol_y) \ar[r] \ar@2{~}[d]
& \Aff(\wt\Fol_x,\wt\Fol_y)
\ar[r] \ar@2{~}[d] &
\Stretch(\wt\Fol_x,\wt\Fol_y) \ar[r] \ar@2{~}[d] & 1 \\
1 \ar[r] & \Isom_h(Q_\gamma) \ar[r] & \Isom(Q_\gamma) \ar[r] & C_\gamma
\ar[r] & 1 
}
$$
where $\Isom_h(Q_\gamma)$ is the subgroup of $\Isom(Q_\gamma)$ preserving
each horizontal leaf, and $C_\gamma = \Isom(Q_\gamma) / \Isom_h(Q_\gamma)$;
the latter group is isomorphic to $D_\infinity$, $\Z$, $\Z/2$, or the
trivial group.

We also have an isomorphism of short exact sequences
$$
\xymatrix{
1 \ar[r] & \Isomplus(\wt\Fol_x,\wt\Fol_y) \ar[r] \ar@2{~}[d]
& \Aff_+(\wt\Fol_x,\wt\Fol_y)
\ar[r] \ar@2{~}[d] &
\Stretch_+(\wt\Fol_x,\wt\Fol_y) \ar[r] \ar@2{~}[d] & 1 \\
1 \ar[r] & \Isom_h(Q_\gamma) \ar[r] & \Isom_+(Q_\gamma) \ar[r] & C_{\gamma+}
\ar[r] & 1 
}
$$
where $\Aff_+(\wt\Fol_x,\wt\Fol_y)$ was defined earlier, 
$\Stretch_+(\wt\Fol_x,\wt\Fol_y)$ is the intersection of
$\Stretch(\wt\Fol_x,\wt\Fol_y) \subset \SL(2,\R)$ with the diagonal
subgroup of $\SL(2,\R)$, $\Isom_+(Q_\gamma)$ is the subgroup of
$\Isom(Q_\gamma)$ preserving the transverse orientation on the horizontal
foliation, and $C_{\gamma+} = \Isom_+(Q_\gamma) / \Isom_h(Q_\gamma)$; in
each of these four cases, the $+$ subscript induces a subgroup of index $\le
2$. The group $\Stretch_+(\wt\Fol_x,\wt\Fol_y)\approx C_{\gamma+}$ is either
$\Z$ or trivial. 

Note that $C_\gamma$ contains the group $\Stab_\T(\gamma) = \{\Phi \in
\Isom(\T) \suchthat \Phi(\gamma)=\gamma\}$, and we shall show in Lemma
\ref{LemmaPeriodic} that this containment has finite index.

\paragraph{Schottky groups and singular \solv-manifolds.} Now consider a
Schottky group $H \subset \MCG$ with limit set $\Lambda
\subset \PMF$, and weak convex hull $\Hull\Lambda$. We have the Cayley
graph $T_H$ and an $H$-equivariant immersion $\rho \from T_H\to\Teich$;
let $\overline T_H = T_H \union \Lambda$. An immediate corollary of
Theorem \ref{theorem:schottkydef} is that $T_H$ and
$\Hull\Lambda$ have finite Hausdorff distance in $\Teich$, and so we
obtain an $H$-equivariant quasi-isometry $\theta \from\Hull\Lambda\to T_H$
which moves points a uniformly bounded distance in $\Teich$. It also
follows from Theorem \ref{theorem:schottkydef} that the map $\theta$
extends, via the identity map on $\Lambda$, to an $H$-equivariant map
$\overline\theta\from \overline\Hull\Lambda \to\overline T_H$ which is
continuous at each point of $\Lambda$. It follows immediately that for
any $\xi\ne\eta \in\Lambda$, the restriction of $\theta$ to the geodesic
$\geodesic{\xi}{\eta}$ under $\theta$ is a quasigeodesic, with
quasigeodesic constants independent of $\xi,\eta$, and the ends of this
quasigeodesic converge in $\overline T_H$ to $\xi,\eta$ respectively;
the unique geodesic in $T_H$ with these endpoints is denoted
$\overline{\xi\eta}$. To summarize:

\begin{itemize}
\item The correspondence between geodesics $\gamma=\geodesic{\xi}{\eta}$
in $\Hull\Lambda$ and geodesics $w=\overline{\xi\eta}$ in $T_H$ is a
bijection. We denote this correspondence by $w=w_\gamma$,
$\gamma=\gamma_w$. Corresponding geodesics $\gamma,w_\gamma$ are uniformly
Hausdorff close in $\Teich$.
\end{itemize}

We need to be a bit more precise. By lifting each geodesic in
$\Hull\Lambda$ to the unit tangent bundle $T^1\T$ we get a closed subset
of
$T^1\T$ invariant under the geodesic flow denoted $T\Hull\Lambda$. The
space $T\Hull\Lambda$ is locally homeomorphic to
$\Lambda\cross\Lambda\cross\R$, i.e.\ it is locally a
Cantor set crossed with the line. We will often confuse a geodesic
$\gamma$ in $\Hull\Lambda$ with its lift to
$T\Hull\Lambda$; these geodesics form a lamination of $T\Hull\Lambda$. As
we have said, the map $\theta\from\Hull\Lambda\to T_H$ lifts to a map
from $T\Hull\Lambda$ to $T_H$, taking each geodesic
$\gamma$ to a quasigeodesic in $T_H$ uniformly Hausdorff close to
$w_\gamma$. But then, by moving the map $T\Hull\Lambda \to T_H$ a bounded
amount, we obtain a continuous map $\Theta\from T\Hull\Lambda \to \T_H$
with the property that $\Theta(\gamma)=w_\gamma$ for each $\gamma$. The
restriction to $\gamma$ is denoted $\Theta_\gamma \from \gamma \to
w_\gamma$, and this is a quasi-isometry with constants independent of
$\gamma$, and we may take $\Theta_\gamma$ to be a homeomorphism.

\begin{proposition}[Comparison of metrics]
\label{proposition:singsol}
Given corresponding geodesics $w$ in $T_H$ and $\gamma=\gamma_w$ in
$T\Hull\Lambda$, there exists a $\pi_1\Sigma$-equivariant, horizontal
respecting quasi-isometry $F_w\from P_w\to Q_{\gamma}$, with
quasi-isometry constants independent of $w$, such that $F_w$ is a lift of
the map $\Theta_\gamma^\inv \from w \to \gamma$.
\end{proposition}

In other words, the natural metric on each hyperplane of the geometric
model space $X_H$ is uniformly quasi-isometric to a singular \solv\
$3$-manifold.  This is the key place where we use the fact that $H$ is a
Schottky group, and not just any free group of pseudo-Anosovs.

\begin{proof} Recall that we have a $\Gamma_H$-equivariant $D$-bundle
$X_H \to T_H$, carrying a $\Gamma_H$-equivariant piecewise Riemannian
metric whose restriction to each fiber is isometric to $\hyp^2$, with
$\Gamma_H$ acting cocompactly. Each of the spaces $P_w$ is embedded in
$X_H$ as the inverse image of $w \subset T_H$.

We now define a $\Gamma_H$-equivariant $D$-bundle $\Xi \to T\Hull\Lambda$
in which each of the singular \solv-manifolds $Q_\gamma$ sits, as follows.

First, note that each point of $T^1\Teich$ corresponds to (the isotopy
class of) a normalized $xy$-structure $(\Fol_x,\Fol_y)$ on $\Sigma$. We
may assemble these structures into an $\MCG$-equivariant $\Sigma$-bundle
over $T^1\Teich$, each fiber equipped with an $xy$-structure in the
appropriate isotopy class, so that the $xy$-strucures vary continuously
as the base point in $T^1\Teich$ varies. Restricting to $T\Hull\Lambda$
we obtain an $H$-equivariant $xy$-$\Sigma$-bundle $\Upsilon \to
T\Hull\Lambda$. 

The universal cover of the $\Sigma$-bundle over $T^1\Teich$ is a
$D$-bundle over $\Teich$, with smoothly varying $xy$-structures on
the fibers, on which the $\pi_1\Sigma$ extension of $\MCG$ acts, namely
the once-punctured mapping class group $\MCG(\Sigma,p)$. Restricting to
$T\Hull\Lambda$ we obtain a $\Gamma_H$-equivariant $xy$-$D$-bundle $\Xi
\to T\Hull\Lambda$. Restricting to any geodesic $\gamma \subset
T\Hull\Lambda$ we obtain the fibration $Q_\gamma \to \gamma$. Note that
the foliation of
$T\Hull\Lambda$ by geodesics lifts to a foliation of $\Xi$ by
\nb{3}manifolds: the \nb{3}manifold over the geodesic $\gamma \subset
T\Hull\Lambda$ is $Q_\gamma$.  The fiberwise $xy$-structures vary
continuously in $\Xi$, and the arc length parameter on geodesics of
$T\Hull\Lambda$ varies continuously; as noted earlier, these data
determine singular \solv-metrics on each $Q_\gamma$, and these metrics
vary continuously in~$\Xi$.

Now lift the $H$-equivariant map $\Theta \from T\Hull\Lambda \to T_H$ to
a $\Gamma_H$-equivariant continuous map of $D$-bundles, $\wt\Theta \from
\Xi\to X_H$, taking each fiber of $\Xi$ homeomorphically to the
corresponding fiber of $X_H$, and taking each singular \solv\ manifold
$Q_\gamma$ to the corresponding singular Riemannian manifold
$P_{w_\gamma}$. 

By cocompactness of the $\Gamma_H$ actions on $\Xi$ and
$X_H$, and by continuity of the $xy$-structures on fibers of $\Xi$ and
the $\hyp^2$ structures on fibers of $X_H$, it follows that $\wt\Theta$
induces quasi-isometries from $\Xi$ fibers to $X_H$ fibers with uniform
quasi-isometry constants; in particular, we get uniform quasi-isometries
from the horizontal sets of $Q_\gamma$ to the horizontal sets of
$P_{w_\gamma}$, over all geodesics $\gamma$ in $T\Hull\Lambda$. Moreover,
since the map $\Theta \from T\Hull\Lambda \to T_H$ is uniformly
quasi-isometric from a geodesic $\gamma$ in $T\Hull\Lambda$ to the
corresponding geodesic $w_\gamma$ in $T_H$, it follows that
$\wt\Theta$ induces a uniform quasi-isometry from the $\abs{dt}$ term in
the metric on $Q_\gamma$ to the $d\tau$ term in the metric on
$P_{w_\gamma}$. This implies that the family of maps $Q_\gamma \to
P_{w_\gamma}$ is uniformly quasi-isometric.
\end{proof}


\section{Quasi-isometries remember the dynamics}
\label{section:dynamics}

Let $H \subset \MCG$ be a Schottky subgroup and let
$\phi \from \Gamma_H \to \Gamma_H$ be a quasi-isometry.  
By Propositions \ref{proposition:hyperplanes:respected},
\ref{proposition:horizontal:respected}, and \ref{proposition:singsol},
to each bi-infinite path $w \in T_H$ there
corresponds a bi-infinite path $w' \in T_H'$, and a horizontal
respecting quasi-isometry $Q_{\gamma_w} \to Q_{\gamma_{w'}}$ of singular
\solv\ spaces. In this section we study such quasi-isometries, and show
that they must coarsely respect certain dynamically defined foliations.

\subsection{The vertical flow on $Q_\gamma$}

Let $\gamma$ be a fixed \Teichmuller\ geodesic in $\Teich(\Sigma)$, and
let $Q_\gamma$ be the associated singular \solv\ metric on $D \cross
\gamma \homeo D \cross \R$. There is a natural flow on $Q_\gamma$ given
by 
$$\Psi_s(x,t)=(x,t+s)$$

This flow is called the \emph{vertical flow}, and its orbits are called 
\emph{vertical geodesics} or \emph{vertical flow lines} in $Q_\gamma$.  
This flow is a pseudo-Anosov flow in the sense
of \cite{FenleyMosher:qfl}. The stable  and unstable foliations of $\Psi$
expand and contract at the uniform rate $e^t$. Singular orbits of $\Psi$
are the same as singular vertical geodesics of $Q_\gamma$.

There are three naturally defined foliations of $Q_\gamma$ which are
invariant under the vertical flow: the codimension-2 foliation by
vertical flow lines; the codimension-1 weak stable foliation; and the
codimension-1 weak unstable foliation. Given a vertical flow line $\ell$,
the {\em weak stable} (resp.\ {\em unstable}) 
leaf containing $\ell$ is the union of
vertical flow lines that asymptote to $\ell$ as $t \to \infinity$ (resp.\
$t\to-\infinity$). Each of the weak stable and unstable foliations has
some nonmanifold leaves, one such leaf containing each singular vertical
geodesic. Note that for each leaf $L$ of the weak stable or unstable
foliation, if $L$ is nonsingular then $L$ is isometric to a hyperbolic
plane, whereas if $L$ is the singular leaf through a singular
vertical geodesic $\ell$ then $L$ is isometric to a union of hyperbolic
half-planes meeting along their common boundary $\ell$; in either case,
$L$ can be expressed as a finite union of hyperbolic planes.

Our goal now is to use pseudo-Anosov dynamics to show that a
horizontal-respecting quasi-isometry between two hyperplanes must also
coarsely respect all of these foliations. In short: the quasi-isometry
remembers the dynamics.

\subsection{Coarse intersection}

We need a basic notion from coarse topology.

\begin{Definition}[Coarse intersection]
A subset $W$ of a metric space $X$ is a \emph{coarse intersection}
of subsets $U,V \subset X$, denoted $W=U\coarsecap V$, if there
exists $C_0$ such that for every $C\geq C_0$ there exists $A=A(C)\geq
0$ so that 
$$d_\Haus(\nbhd_C(U)\cap \nbhd_C(V),W)\leq A
$$  
Note that although such a set $W$ may not exist, when it does
exist then any two such sets are a bounded Hausdorff distance from
each other. The function $A(C), C \ge C_0$ is called the \emph{coarse
intersection function}.
\end{Definition}

We will need the following fact, which is an elementary consequence of the 
definitions.

\begin{lemma}
\label{lemma:intersection2}
For any quasi-isometry $f\from X\to Y$ of metric spaces, and $U,V \subset 
X$, if $U \coarsecap V$ exists then $f(U \coarsecap V)$ is a coarse 
intersection of $f(U)$, $f(V)$, with coarse intersection function depending
only on the quasi-isometry constants for $f$ and the coarse intersection
function for $U$ and $V$.
\qed\end{lemma}

\subsection{The dynamically defined foliations are respected}

In this subsection we prove the key proposition: 

\begin{proposition}[Stable and unstable foliations respected]
\label{proposition:foliations:respected} 
Let $\gamma, \gamma'$ be cobounded \Teichmuller\ geodesics. Then any
horizontal-respecting quasi-isometry
$\phi\from Q_\gamma\to Q_{\gamma'}$ coarsely respects the stable and
unstable foliations of the vertical flows on $Q_\gamma,Q_{\gamma'}$.
Also, $\phi$ coarsely respects the patterns of singular stable and
unstable leaves.
\end{proposition}

Every vertical flow line can be realized as the coarse intersection of a
stable and unstable leaf, with uniform coarse intersection function 
independent of the choice of flow line. Proposition
\ref{proposition:foliations:respected} and Lemma
\ref{lemma:intersection2} therefore imply that $\phi$ coarsely respects
the  collection of vertical geodesics. A similar argument works for the
collection of singular geodesics, using singular stable and unstable
leaves.  We record this as: 

\begin{corollary}[Vertical flow lines respected]
\label{corollary:vertical:respected}
The quasi-isometry $\phi$ of Proposition
\ref{proposition:foliations:respected} coarsely 
respects the foliations of vertical flow lines in $Q_{\gamma}$ and
$Q_{\gamma'}$, as well as the patterns of singular vertical lines.
\end{corollary}

\begin{proof}[Proof of Proposition \ref{proposition:foliations:respected}] 
The horizontal foliation of $Q_\gamma$ is an example of a
\emph{uniform foliation}, which means that any two leaves have finite
Hausdorff distance. Any map between two horizontal leaves which moves
points a bounded distance in $Q_\gamma$ is a quasi-isometry between those
leaves. It follows that there is a canonical coarse equivalence class of
quasi-isometries between any two leaves, and moreover the composition of
two such quasi-isometries is another one. Each leaf is quasi-isometric to
the hyperbolic plane $\hyp^2$ and its Gromov boundary is a
circle, so there is a canonical identification of all the circles at
infinity to a single circle which we denote $SQ_\gamma$. These facts were
noted by Thurston in \cite{Thurston:FoliationsAndCircles}.

It is well-known that a 
quasi-isometry $\hyp^2\to \hyp^2$ induces a homeomorphism $\partial
\hyp^2 \to \partial \hyp^2$ between the circles at infinity, and that
coarsely equivalent quasi-isometries induce the same boundary map.  
Any horizontal respecting quasi-isometry $Q_\gamma \to Q_{\gamma'}$
therefore induces a homeomorphism $SQ_\gamma \to SQ_{\gamma'}$ between
their respective circles at infinity. The underlying idea of our proof is
to find additional quasi-isometrically invariant structures on $SQ_\gamma$
and
$SQ_{\gamma'}$ which encode the stable and unstable foliations, and use
this information to prove quasi-isometric invariance of these
foliations.

To proceed with the proof we need some notation. Let $D_t$ denote the
plane $D \cross t \subset D \cross \R = Q_\gamma$. The plane $D_t$ comes
equipped with an $xy$-structure, as explained in
\S\ref{section:singularSolv}. For each $s,t$ let $\phi_{st} \from D_s \to
D_t$ be the map $(p,s) \to (p,t)$, $p \in D$. In other words, $\phi_{st}$
flows along vertical flow segments of $Q_\gamma$ from $D_s$ to $D_t$;
each such flow segment has length $\abs{s-t}$, and so $\phi_{st}$ is a
quasi-isometry in the canonical coarse equivalence class from $D_s$ to
$D_t$ as discussed above. In fact, the map $\phi_{st}$ is an $xy$-affine
homeomorphism with stretch factor $e^{\abs{s-t}}$, implying that
$\phi_{st}$ is $e^{\abs{s-t}}$-bilipschitz. Note that $\phi_{tu}
\composed \phi_{st} = \phi_{su}$, for all $s,t,u \in \R$. 

The $xy$-metric on $D_t$ is a $\CAT(0)$ metric, and it is also Gromov
hyperbolic because $D_t$ is quasi-isometric to $\hyp^2$. The boundary
$\bdy D_t$ is a circle, and each quasi-isometry $\phi_{st}\from D_s\to D_t$
induces a homeomorphism $\partial \phi_{st}\from \partial D_s\to \partial
D_t$.  

For any two points $\xi,\eta \in \bdy D_t$ there exists an
$xy$-geodesic with endpoints $\xi,\eta$, and any two such geodesics are
the boundary of an isometric embedding of $\R
\cross [a,b]$ for some $[a,b] \in \R$ (see \cite{BridsonHaefliger}). Since
$D_t$ is not the Euclidean plane but has a cocompact isometry group,
there is an upper bound on $b-a$ independent of $\xi,\eta$. This bound is
moreover independent of $t$, because coboundedness of the \Teichmuller\
geodesic $\gamma$ and compactness of the fibers of $T^1\T$ together imply
that the collection of locally $\CAT(0)$ metric spaces $D_t /
\pi_1(\Sigma)$ lies in a compact space of locally $\CAT(0)$ metrics
(this is the one place in the proof where we use coboundedness of
the \Teichmuller\ geodesics in the weak convex hull of $H$).  
Let $R_0$ be a $t$-independent bound for $\abs{b-a}$.

For each bi-infinite geodesic $\theta$ on $D_0$, the
set 
$$V = V_\theta = \Union_{t \in \R}\{\phi_{0t}(\theta)\}
$$
is called a \emph{vertical plane} in $Q_\gamma$. Note that if $\bdy\theta
= \{\xi,\eta\} \subset \bdy D_0$ then $\bdy\phi_{0t}(\theta) =
\{\bdy\phi_{0t}(\xi),\bdy\phi_{0t}(\eta)\} \subset \bdy D_t$, 
and so to $V$ there
is associated a unique pair of points in the circle $SQ_\gamma$ which we
call the \emph{endpoints} of $V$. Moreover, the discussion in the
previous paragraph shows that the Hausdorff distance between any two
vertical planes of $Q_\gamma$ with the same endpoints $\xi,\eta$ is at
most $R_0$.

The first structure which $\phi\from Q_\gamma\to Q_{\gamma'}$ must coarsely 
respect is the collection of vertical planes.

\begin{claim}[Vertical planes preserved]
\label{claim:planes1}
Let $\phi \from Q_\gamma \to Q_{\gamma'}$ be a horizontal-respecting
quasi-isometry.  If $V$ is any vertical plane in $Q_\gamma$, then 
$\phi(V)$ is a bounded Hausdorff distance from some vertical plane 
$V'$ in $Q_{\gamma'}$.
\end{claim}

To prove Claim \ref{claim:planes1}, let $t \to t'$ be a bijective
quasi-isometry of $\R$ such that $\phi(D_t)$ is Hausdorff close to
$D'_{t'}$, with a uniform Hausdorff constant; we assume that the
parameterizations are chosen so that $0 \to 0$ under this
quasi-isometry. Composing the map $D_t \to \phi(D_t)$ with a uniformly
finite distance map $\phi(D_t) \to D'_{t'}$, we obtain a quasi-isometry
$\psi_t \from D_t\to D'_{t'}$ whose quasi-isometry constants are
independent of $t$.  

Fix a vertical plane $V=V_\theta$ and for each $s$ consider the geodesic
$\theta_s = V \intersect D_s$. Its image quasigeodesic $\psi_s(\theta_s)$
is uniformly Hausdorff close in $D'_{s'}$ to some geodesic
$\theta'_{s'}$, in particular $\phi_0(\theta_0)$ is close to $\theta'_0$.
We need only show that for each $s$, $\theta'_{s'}$ is uniformly
Hausdorff close to
$\phi'_{0s'}(\theta'_0)$, for then we can set $$V' =
\union_{s'}\{\phi'_{0s'}(\theta'_0)\}$$ and it follows that $\phi(V)$ is
uniformly Hausdorff close to $V'$.

Since any two geodesics in $D'_{s'}$ which are Hausdorff close are
$R_0$-Hausdorff close, we need only show that the Hausdorff distance
between $\theta'_{s'}$ and $\phi'_{0s'}(\theta'_0)$ is finite, in other
words these two geodesics have the same endpoints in $\bdy D'_{s'}$. But
this is an immediate consequence of the coarse commutativity of the
following diagram:
$$\xymatrix{
D_0 \ar[r]^{\psi_0} \ar[d]_{\phi_{0s}} & D'_{0} \ar[d]^{\phi'_{0s'}} \\
D_s \ar[r]_{\psi_s} & D'_{s'}
}
$$

This finishes the proof of Claim \ref{claim:planes1}.  We now find
finer structures which must be coarsely respected by $\phi$.
\bigskip

Each vertical plane $V$ comes equipped with a \emph{horizontal
foliation}, obtained by intersecting the plane with the horizontal
foliation $\{D_t\}$ of $Q_\gamma$. Denote the leaves of this foliation
by 
$$\theta_t= D_t\intersect V, \quad t \in \R$$ 
so that $\theta_t =
\phi_{st}(\theta_s)$ for any $s,t \in \R$. The horizontal foliation
$\{\theta_t\}$ has a transverse orientation, pointing in the direction
of increasing $t$. Note that the Hausdorff distance in $V$ between
$\theta_t$ and $\theta_s$ is exactly $\abs{t-s}$. There is a projection
$\pi \from V \to \R$ with $\theta_t =
\pi^\inv(t)$. Define a \emph{quasivertical line} in $V$ to be a subset $L
\subset V$ such that the projection $\pi \from L \to \R$ is a
quasi-isometry.

We divide the collection of vertical planes into three types: stable,
unstable, and doubly unstable. A \emph{stable} vertical plane is one which
is contained in a leaf of the stable foliation on $Q_\gamma$; thus it is
either a regular leaf of the stable foliation, or it is a union of two
half-planes of a singular leaf. Similarly an \emph{unstable} vertical
plane is one contained in an unstable leaf. All other vertical planes 
are called \emph{doubly unstable} vertical planes.

The three types of vertical planes---stable, unstable, and
doubly unstable---can be distinguished from each other by horizontal
respecting quasi-isometries which respect the transverse orientation, by
observing the asymptotic behavior of quasivertical lines. To make this
more precise,  say that two quasivertical lines $L,L' \subset P$ are
\emph{upward Hausdorff close} if the Hausdorff distance between $L
\intersect
\pi^\inv[0,\infinity)$ and $L' \intersect \pi^\inv[0,\infinity)$ is
finite; \emph{downward Hausdorff close} is similarly defined using
$(-\infinity,0]$. 

\begin{claim}
\label{claim:planes2}
Let $V$ be a vertical plane. Then:
\begin{description}
\item[Stable:] If $V$ is a stable plane, 
then any two quasivertical lines
in $V$ are upward Hausdorff close but not downward Hausdorff close.
\item[Unstable:] If $V$ is an unstable plane, then any two
quasivertical lines in $V$ are downward Hausdorff close but not upward
Hausdorff close.
\item[Doubly unstable:] If $V$ is a doubly unstable plane then there exist
two quasivertical lines in $V$ which are neither upward Hausdorff close
nor downward Hausdorff close.
\end{description}
\end{claim}

To prove Claim \ref{claim:planes2} when $V$ is stable or unstable, 
observe first that $V$
with its horizontal foliation $\{\theta_t\}$ is isometric to the
hyperbolic plane $\hyp^2$ in the upper half plane model, with the
``horizontal'' horocyclic foliation centered on the point $\infinity$. If
$V$ is stable (resp.\ unstable) then the transverse orientation points
towards $\infinity$ (resp.\ away). Next observe that any quasivertical
line in $\hyp^2$ is a quasigeodesic with one endpoint at $\infinity$, and
so any two quasivertical lines are Hausdorff close in the direction of
$\infinity$. 

To prove Claim \ref{claim:planes2} when $V$ is doubly unstable, observe
that $\theta_0=V\intersect D_0$ is an $xy$-geodesic in $D_0$ which is not
contained in a leaf of either the $x$-foliation or the $y$-foliation on
$D_0$. There exists, therefore, two points $p,q \in \theta_0$ which do not
lie in the same leaf of either the $x$-foliation or the $y$-foliation on
$D_0$. The vertical flow lines $p \cdot \R$, $q \cdot \R$ in $Q_\gamma$
are evidently neither upward nor downward Hausdorff close in the singular
\solv-metric on $Q_\gamma$, and so the same is true in $V$.  This proves 
Claim \ref{claim:planes2}.

The same discussion holds, of course, in $Q_{\gamma'}$. Now consider a
horizontal-respecting quasi-isometry $\phi \from Q_\gamma \to
Q_{\gamma'}$. By reversing the upward orientation in $Q_{\gamma'}$ if
necessary, we may assume that $\phi$ respects the upward
orientation. Let $V$ be any vertical plane in $Q_\gamma$. We have shown
in Claim \ref{claim:planes1} that $\phi(V)$ is Hausdorff close to some
vertical plane $V'$ in
$Q_{\gamma'}$. Composing the map $V\to\phi(V)$ with any finite distance
map $\phi(V) \to V'$, we obtain therefore a quasi-isometry $\psi\from V \to
V'$. Each horizontal leaf of $V$ (resp.\ $V'$) is a coarse intersection
of $V$ with a horizontal leaf of $Q_\gamma$ (resp.\ $Q_{\gamma'}$), and
it follows that $\psi$ coarsely respects the horizontal foliations and
their transverse orientations.

Claim \ref{claim:planes2} shows manifestly that stable vertical planes
are coarsely respected by $\phi \from Q_\gamma \to Q_{\gamma'}$, and
similarly for unstable vertical planes. To finish the proof of
Proposition \ref{proposition:foliations:respected}, given two stable
(resp.\ unstable) vertical planes $V_1,V_2$ in $Q_\gamma$ or in
$Q_{\gamma'}$, we must give a quasi-isometrically invariant property which
characterizes $V_1,V_2$ lying in the same weak stable
(resp.\ unstable) leaf. Namely,
$V_1,V_2$ lie in the same leaf if and only if, for any $t \in \R$, the
triple coarse intersection $V_1\coarsecap V_2 \coarsecap D_t$ is
unbounded.
\end{proof}

\paragraph{Remark.} There is an alternative ``dynamical'' proof of
Proposition \ref{proposition:foliations:respected} which we worked out
first, and which parallels the proof of the analogous proposition in
\cite{FarbMosher:ABC}. The latter is concerned with the solvable Lie
groups associated to geodesics in the symmetric space of $\GL(n,\R)$, and
the natural Anosov flows on these solvable Lie groups; the proof in
\cite{FarbMosher:ABC} is an easy consequence of the shadowing lemma for
Anosov flows. Unfortunately the shadowing lemma is false in pseudo-Anosov
dynamics, complicating the situation drastically. There is an alternative
shadowing theory for pseudo-Anosov dynamical systems, developed in
\cite{Handel:globalshadowing}, \cite{Handel:entropy}, and
\cite{Mosher:Equivariant}, and this was used in our original proof of
Proposition \ref{proposition:foliations:respected}. The proof given above
for Proposition \ref{proposition:foliations:respected} entirely avoids
these issues by using Gromov hyperbolicity of the horizontal leaves---a
fact which was not available in \cite{FarbMosher:ABC}, where the
horizontal leaves are Euclidean.


\section{Periodic hyperplanes}
\label{section:periodic}

A hyperplane $P_w$ in $X_{H}$ is a \emph{periodic hyperplane} if $w$ is  a
periodic geodesic in $T_H$, that is, the subgroup $C_w$ of $H$ stabilizing
$w$ is infinite cyclic. Periodic hyperplanes $P_w$ are special in
that they admit \emph{a priori} extra isometries, above and beyond the
isometric action of $\pi_1(\Sigma)$: the subgroup
$\pi_1(\Sigma)\semidirect C_w$ of $\Gamma_H$ acts isometrically on $P_w$. 

Our goal now is to show that horizontal-respecting quasi-isometries of
periodic hyperplanes must remember the extra symmetries.

Recall that $w \leftrightarrow \gamma_w$ is a bijection between geodesics
in the tree $T_H$ and \Teichmuller\ geodesics in $\Hull\Lambda$. 

\begin{proposition}[Periodic hyperplanes preserved]
\label{proposition:periodicity:preserved} 
Given $w,w'$ geodesics in $T_H$, suppose that there exists a horizontal
respecting quasi-isometry between hyperplanes
$\phi\from P_w \to P_{w'}$. Then:
\begin{enumerate}
\item $w$ is periodic if and only if $w'$ is periodic. 
\item If $w,w'$ are periodic, and if $\gamma, \gamma' \subset
\Hull\Lambda$ are the geodesics in $\T$ associated to $w,w'$ respectively,
then $\gamma, \gamma'$ are periodic and there is a singular \solv-isometry
$\Phi \from Q_\gamma \to Q_{\gamma'}$ such that the following diagram
coarsely commutes:
$$
\xymatrix{
P_w \ar[r]^\phi \ar[d]_{F_w}      &  P_{w'} \ar[d]^{F_{w'}} \\
Q_\gamma \ar[r]^\Phi &  Q_{\gamma'}
}
$$
That is, $d(F_{w'}(\phi(x)),\Phi(F_w(x))) \le A$ where $A$ depends only
on the quasi-isometry constants of $\phi$, not on $w$ or~$w'$.
\end{enumerate}
\end{proposition}

We start by reducing the proposition to a statement about singular
\solv-manifolds by using the following result, which shows that various
competing notions of ``periodicity'' are in fact equivalent:

\begin{lemma}[Characterizing periodicity]
\label{LemmaPeriodic}
Given a $T_H$ geodesic $w$ and the corresponding \Teichmuller\ geodesic
$\gamma=\gamma_w$, the following are equivalent:
\begin{enumerate}
\item $P_w$ is periodic, meaning that the group $C_w$ is infinite.
\item $\gamma$ is periodic in $\T$, meaning that the group
$\Stab_\T(\gamma)$ is infinite.
\item $Q_{\gamma}$ is periodic, meaning that the group $C_{\gamma} =
\Isom(Q_\gamma) / \Isom_h(Q_\gamma)$ is infinite.
\end{enumerate}
\end{lemma}

\begin{proof} We clearly have inclusions $C_w \subset \Stab_\T(\gamma_w)
\subset C_\gamma$ and so (1) implies (2) implies (3).

To prove (2) implies (1), suppose that $\gamma$ is periodic in
$\T$, that is, $\gamma$ is the axis of some pseudo-Anosov element $\phi \in
\MCG$. Choose a point $p \in\gamma$, and consider the sequence $\phi^n(p)$,
$n>0$. Since $H$ acts cocompactly on $\Hull\Lambda$, there is a sequence
$\psi_n
\in H$ such that $\{\psi_n \composed \phi^n(p) \suchthat n>0\}$ is a
bounded subset of $\T$. Since $\MCG$ acts properly on $\T$, there exist $m >
n > 0$ such that $\psi_n \composed \phi^n = \psi_m \composed \phi^m$. It
follows that $\phi^{m-n} = \psi_m^\inv \psi_n^{\vphantom{\inv}} \in H$, and
so $w$ is periodic in $T_H$.

To prove (3) implies (2), choose a horizontal leaf $D_\gamma$ of $Q_\gamma$
and so we have the split exact sequence
$$1 \to \Isomplus(D_\gamma) \approx \Isom_h(Q_\gamma) \to \Isom_+(Q_\gamma)
\to C_{\gamma+} \to 1
$$
where $C_{\gamma+}$ is an infinite cyclic subgroup of index $\le 2$ in
$C_\gamma$. We also have a finite index inclusion $\pi_1(\Sigma)
\subgroup \Isomplus(D_\gamma)$. Since $C_{\gamma+}$ acts by automorphisms of
$\Isomplus(D_\gamma)$ it follows that $C_{\gamma+}$ has a finite index
subgroup stabilizing $\pi_1(\Sigma)$; this subgroup is clearly
identified with a subgroup of $\Stab_\T(\gamma_w)$ and so the latter is
infinite.
\end{proof}

Note from the proof that each of the inclusions
$C_w \subset \Stab_\T(\gamma_w) \subset C_\gamma$ has finite index;
examples may be constructed in which any number of these inclusions is
proper.

By combining Lemma \ref{LemmaPeriodic} with Proposition
\ref{proposition:singsol}, it follows immediately that Proposition
\ref{proposition:periodicity:preserved} is reduced to the following:

\begin{lemma}
\label{LemmaTeichPeriodic}
Given geodesics $\gamma,\gamma'$ in $\T$, if $\gamma$ is
periodic, and if there exists a horizontal respecting quasi-isometry 
$Q_\gamma \to Q_{\gamma'}$, then $\gamma'$ is periodic and there exists a
singular \solv\ isometry $Q_\gamma \to Q_{\gamma'}$.
\end{lemma}

\begin{Remark}
Our proof of Lemma \ref{LemmaTeichPeriodic} uses both Thurston's 
hyperbolization theorem for mapping tori of pseudo-Anosov homeomorphisms
\cite{Otal:fibered} as well as the geodesic pattern rigidity theorem of
R.\ Schwartz \cite{Schwartz:Symmetric}, and to apply these results we
need to invoke periodicity of $\gamma$. Existence of a horizontal
respecting quasi-isometry between $Q_\gamma$ and $Q_{\gamma'}$ ought to
imply existence of a singular \solv-isometry, without assuming
periodicity of $\gamma$; however, we do not know a proof.
\end{Remark}

\subsection{Rigidity of periodic hyperplanes}

The main step in the proof of Lemma \ref{LemmaTeichPeriodic} is
Proposition \ref{prop:periodic:rigidity} below, which shows that a periodic
singular
\solv-manifold $Q_\gamma$ has the rigidity 
property that its horizontal respecting
quasi-isometry group equals its isometry group. Here is the basic setup.

For each periodic geodesic $\gamma$ in $\Teich$, fix a base horizontal
leaf $D_\gamma \subset Q_\gamma$.  Recall that $\pi_1(\Sigma)$ acts as a
deck transformation group on $D_\gamma$. The singular \solv-metric on
$Q_\gamma$ induces an $xy$-structure on $D_\gamma$ with respect to which
we have an inclusion $\pi_1(\Sigma) \hookrightarrow
\Isomplus(D_\gamma)$. There is a commutative diagram of short exact
sequences $$
\xymatrix{
1 \ar[r]  & \Isomplus(D_\gamma) \ar[r] & \Isom(Q_\gamma) \ar[r] & C_\gamma
\ar[r] & 1 \\
1 \ar[r]  & \pi_1(\Sigma) \ar[r] \ar@{^{(}->}[u] &\Gamma_\gamma \ar[r]
\ar@{^{(}->}[u] & \Stab_\T(\gamma) \ar[r] \ar@{^{(}->}[u] & 1
}
$$
where each vertical arrow is a finite index inclusion. The group $C_\gamma$
is either $\Z$ or $D_\infty$. The quotient $Q_\gamma /
\Isom(Q_\gamma)$ is a closed 3-dimensional orbifold, equipped with an
orbifold fibration over the circle (when $C_\gamma=\Z$) or the interval
orbifold (when $C_\gamma = D_\infty$); in the latter case we choose
$D_\gamma$ so that it projects to a generic point of the interval orbifold,
i.e.\ so that the stabilizer of $D_\gamma$ in $\Isom(Q_\gamma)$ equals
$\Isomplus(D_\gamma)$. It follows that the generic fiber of the quotient
3-orbifold $Q_\gamma /
\Isom(Q_\gamma)$ is the closed 2-orbifold $\O_\gamma = D_\gamma /
\Isomplus(D_\gamma)$ equipped with the quotient $xy$-structure. 

Note that the quotient 3-orbifold $Q_\gamma / \Gamma_\gamma$ also
fibers over the circle or the interval orbifold, with fiber $\Sigma$. There
is a finite covering map $\Sigma \to \O_\gamma$ and the fiber monodromy map
on $\O_\gamma$ lifts to a pseudo-Anosov homeomorphism $\phi \from \Sigma
\to \Sigma$.

Let $\QI_h(Q_\gamma)$ be the subgroup of $\QI(Q_\gamma)$
represented by horizontal respecting quasi-isometries of
$Q_\gamma$.  Clearly there is a homomorphism
$$i \from \Isom(Q_\gamma) \to \QI_h(Q_\gamma)
$$

\begin{proposition}[Rigidity of periodic hyperplanes]
\label{prop:periodic:rigidity}
If $\gamma$ is periodic then the homomorphism $i \from
\Isom(Q_\gamma)\to \QI_h(Q_\gamma)$ is an isomorphism.
\end{proposition}

\begin{proof}
A nontrivial isometry $\phi$ of $Q_\gamma$ must map some vertical
geodesic to a different vertical geodesic.  Since any two vertical
geodesics in $Q_\gamma$ have infinite Hausdorff distance, $\phi$ is an
infinite distance from the identity, so that $i$ is injective
(injectivity of $i$ therefore is true regardless of periodicity).

We now prove that $i$ is surjective: every horizontal respecting
quasi-isometry of $Q_\gamma$ is a bounded distance from an isometry.

Thurston's geometrization theorem for pseudo-Anosov mapping tori gives a
hyperbolic structure on the 3-dimensional orbifold $Q_\gamma
/\Isom(Q_\gamma)$. This yields a properly discontinuous, cocompact,
isometric, faithful action 
$$h \from \Isom(Q_\gamma) \to\Isom(\hyp^3)
$$ 
and a quasi-isometric, $h$-equivariant
homeomorphism 
$$q \from Q_\gamma \to \hyp^3
$$

For each singular vertical geodesic $\ell \subset Q_\gamma$, the image
$q(\ell)$ is a quasigeodesic in $\hyp^3$, which by the Morse-Mostow Lemma
is a bounded  Hausdorff distance (not depending on $\ell$) from a unique
geodesic in $\hyp^3$; let $L$ denote the set of all such geodesics in
$\hyp^3$. Let $\Isom(\hyp^3,L)$ be the group of isometries of $\hyp^3$
which permute the collection $L$. Also, let $F$ be the foliation of
$\hyp^3$ obtained by pushing forward via $q$ the horizontal foliation of
$Q_\gamma$, and let $\Isom(\hyp^3,L,F)$ be the subgroup of
$\Isom(\hyp^3,L)$ which coarsely respects $F$. Note that
$\image(h) \subset \Isom(\hyp^3,L,F)$. Now we show that $h$ factors
through~$i$.

For each horizontal respecting quasi-isometry $\psi$ of $Q_\gamma$, we obtain
a quasi-isometry $q\psi q^\inv$ of $\hyp^3$, inducing a homomorphism
$\hat q \from \QI_h(Q_\gamma) \to \QI(\hyp^3)$. Since
$\psi$ coarsely respects singular vertical geodesics and horizontal
leaves in $Q_\gamma$, it follows that $q \psi q^\inv$ coarsely respects
$L$ and $F$. 

Since $q \psi q^\inv$ coarsely respects $L$, since $L$ is invariant under
the cocompact isometry group $\image(h)$, and since there are only
finitely many orbits of the action of $\image(h)$ on $L$, we may
directly apply the main theorem of \cite{Schwartz:Symmetric} which says in
this setting that $q \psi q^\inv$ is a bounded distance (not uniformly so)
from a unique isometry of $\hyp^3$ which strictly respects $L$. 

It follows that the image of $\hat q$ is contained in $\Isom(\hyp^3,L)$,
and in fact in $\Isom(\hyp^3,L,F)$. It's evident that the composition
$$\Isom(Q_\gamma) \xrightarrow{i} \QI_h(Q_\gamma) \xrightarrow{\hat q}
\Isom(\hyp^3,L,F) \subset \Isom(\hyp^3)
$$
is identical with the homomorphism $h$. The homomorphism $\hat q$ is
obviously injective, and so to prove surjectivity of $i$ it suffices to
show that $\image(h) = \Isom(\hyp^3,L,F)$.

Consider the short exact sequence
$$1 \to \Isomplus(D_\gamma) \to \Isom(Q_\gamma) \to C_\gamma \to 1
$$
Since $\Isomplus(D_\gamma)$ is normal in $\Isom(Q_\gamma)$, and
since $\Isom(Q_\gamma)$ is identified via $h$ with a finite index subgroup
of $\Isom(\hyp^3,L,F)$ (both being discrete and cocompact on $\hyp^3$), it
follows that
$\Isomplus(D_\gamma)$ has a normalizer of finite index in
$\Isom(\hyp^3,L,F)$; choose coset representatives
$g_1,\ldots,g_n$ of the normalizer. Each leaf of $F$ is coarsely
equivalent to $\Isomplus(D_\gamma)$, and $g_1,\ldots,g_n$ coarsely respects
$F$, and so each conjugate subgroup 
$$g_1 \Isomplus(D_\gamma) g_1^\inv, \ldots, g_n
\Isomplus(D_\gamma) g_n^\inv
$$
is coarsely equivalent to $\Isomplus(D_\gamma)$. 

Now we apply an elementary lemma of \cite{MosherSageevWhyte:QuasiTreeTwo}
which says that for a
finite collection of subgroups in a finitely generated group, the coarse
intersection of those subgroups is coarsely equivalent to their
intersection.  The intersection of the above conjugates of
$\Isomplus(D_\gamma)$ is therefore coarsely equivalent to $D_\gamma$. 

Another elementary lemma of \cite{MosherSageevWhyte:QuasiTreeTwo} says that
in a finitely generated group, given subgroups $A \subset B$, if $A,B$ are
coarsely equivalent, then $A$ has finite index in $B$. It follows that the
intersection of the conjugates of $\Isomplus(D_\gamma)$ in $\Isom(\hyp^3,L,F)$
has finite index in $\Isomplus(D_\gamma)$. Thus we obtain a normal subgroup $N
\subset\Isom(\hyp^3,L,F)$ of finite index in $\Isomplus(D_\gamma)$. The
quotient group $\Isom(\hyp^3,L,F) / N$ is a finite index supergroup of
$C_\gamma =
\Isom(Q_\gamma) / \Isomplus(D_\gamma)$, and so the quotient is virtually
cyclic. Since $\Isom(\hyp^3,L,F)$ is a \nb{3}orbifold group it follows that
the quotient $\Isom(\hyp^3,L,F)/N$ is either $\Z$ or $D_\infinity$. 

The
orbifold $\hyp^3 / \Isom(\hyp^3,L,F)$ therefore fibers over the
circle or the interval orbifold, with generic fiber $\O$, and with
$\pi_1(\O)$ identified with $N$. Lifting this fibration to $\hyp^3$ we
obtain an
$\Isom(\hyp^3,L,F)$ equivariant fibration coarsely equivalent to
$F$. We may therefore replace $F$ with this fibration, and so $F$ is
strictly invariant under $\Isom(\hyp^3,L,F)$. The monodromy map on
$\O$ is pseudo-Anosov. There is a finite index covering map $Q_\gamma /
\Isom(Q_\gamma) \to \hyp^3 / \Isom(\hyp^3,L,F)$, taking fibration to
fibration. By uniqueness of pseudo-Anosov homeomorphisms in their
isotopy classes \cite{FLP}, it follows that the stable and unstable
measured foliations for the monodromy map of $\O$ lift to the stable and
unstable measured foliations for the monodromy map on the generic fiber
$\O_\gamma$ of
$Q_\gamma / \Isom(Q_\gamma)$. But this shows that
$\Isom(\hyp^3,L,F)$ acts isometrically on $Q_\gamma$, proving that
$\image(h) = \Isom(\hyp^3,L,F)$.

This proves Proposition \ref{prop:periodic:rigidity}.
\end{proof}

\subsection{Proof of Lemma \ref{LemmaTeichPeriodic}}

Assume $Q_\gamma$ is periodic and $\phi \from Q_\gamma \to Q_{\gamma'}$ is a horizontal
respecting quasi-isometry. We'll construct an isometry $\Phi \from Q_\gamma
\to Q_{\gamma'}$, the existence of which implies that $Q_{\gamma'}$ is
periodic. 

Applying Proposition \ref{proposition:foliations:respected}, we may move
$\phi$ a bounded distance so that $\phi$ is a homeomorphism, respecting
the horizontal foliations and the weak stable and unstable foliations. We
may assume that the base horizontal leaves $D_\gamma \subset Q_\gamma$,
$D_{\gamma'} \subset Q_{\gamma'}$ are chosen so that
$\phi(D_\gamma)=D_{\gamma'}$.

Let $A_\gamma \from \pi_1(\Sigma) \to \Isom(Q_\gamma)$ and $A_{\gamma'} \from
\pi_1(\Sigma) \to \Isom(Q_{\gamma'})$ be the standard isometric actions,
preserving each horizontal leaf.
Let $\QIMap_h(Q_{\gamma})$ be the subsemigroup of $\QIMap(Q_\gamma)$ that
coarsely respects the horizontal foliation of $Q_\gamma$, and let
$$B_{\gamma} = \phi^\inv \composed A_{\gamma'} \composed \phi \from \pi_1(\Sigma)
\to \QIMap_h(Q_{\gamma}) 
$$
be the conjugated action (its really an action, not just a quasi-action,
because $\phi$ is a homeomorphism). Note that $B_\gamma$ preserves each
horizontal leaf of $Q_\gamma$ as well as the strong stable and unstable
foliations in that leaf; however, $B_\gamma$ does not a priori preserve the
invariant measures on those foliations.

Applying Proposition \ref{prop:periodic:rigidity}, we conclude
that $B_\gamma$ is a bounded distance from an isometric action, that is,
there exists a homeomorphism $\xi \from Q_\gamma \to Q_\gamma$ which moves each
point a uniformly bounded distance, such that $\xi^\inv \composed B_\gamma
\composed
\xi$ is an isometric action of $\pi_1(\Sigma)$ on $Q_\gamma$. The action
$\xi^\inv \composed B_\gamma \composed  \xi$ also preserves each horizontal
leaf of $Q_\gamma$ and the strong stable and unstable foliations in that leaf.

We claim that $B_\gamma = \xi^\inv\composed B_\gamma \composed \xi$. Consider a
horizontal leaf $L$ of $Q_\gamma$ and a point $x \in L$. Let $\ell^s,
\ell^u$ be the strong stable and unstable leaves in $L$ passing through
$x$, and so $x = \ell^s \intersect \ell^u$. Consider an element $\beta \in
B_\gamma(\pi_1(\Sigma))$. We know that $\beta(L)=\xi^\inv\beta\xi(L)=L$. We
also know that $\beta(\ell^s)$ and $\xi^\inv \beta \xi(\ell^s)$
are both strong stable leaves in $L$, and they are a bounded distance
from each other, implying that $\beta(\ell^s) =
\xi^\inv\beta\xi(\ell^s)$. Similarly,
$\beta(\ell^u)=\xi^\inv\beta\xi(\ell^u)$. Therefore,
$\beta(x)=\xi^\inv\beta\xi(x)$, proving the claim.

It follows from the claim that $B_\gamma$ does, in fact, preserve the
invariant measures on the strong stable and unstable foliations in each
horizontal leaf of $Q_\gamma$.

Fix a horizontal leaf $L \subset Q_\gamma$ and let $L' = \phi(L) \subset
Q_{\gamma'}$.  So $\phi$ takes the strong stable and unstable foliations in $L$
to those in $L'$. Let $f^s, f^u$ be the strong stable and unstable measured
foliations in $L$, and let $f'{}^s, f'{}^u$ the strong stable and unstable
measured foliations in $L'$. The map $\phi$ pushes the transverse measures
on
$f^s,f^u$ forward to new transverse measures on $f'{}^s, f'{}^u$. 

Since $B_\gamma = \phi^\inv A_{\gamma'} \phi$ acts isometrically on the
transverse measures of $f^s,f^u$, it follows that $A_{\gamma'}$ acts
isometrically on both the old and the new transverse measures on $f'{}^s,
f'{}^u$. However, $f'{}^s, f'{}^u$ are \emph{uniquely ergodic} with respect
to the $A_{\gamma'}$ action, i.e.\ they have projectively unique transverse
measures invariant under the $A_{\gamma'}$ action.  This follows from the
fact that the \Teichmuller\ geodesic $\gamma'$ is cobounded, together with
a theorem of H.\ Masur that if $\xi \in \PMF$ is not uniquely ergodic then
any \Teichmuller\ ray with ending foliation $\xi$ is not cobounded
\cite{Masur:UniquelyErgodic}. Hence  the old and new transverse measures on
$f'{}^s, f'{}^u$ differ by multiplicative constants. These multiplicative
constants are inverses to each other, because the action $A_{\gamma'}$ has
cofinite area. It follows that $\phi$ restricts to an $xy$-affine
isomorphism from $D_\gamma$ to $D_{\gamma'}$. This implies in turn that the
quasi-isometry $\phi \from Q_\gamma \to Q_{\gamma'}$ may be altered a
bounded amount, moving each point a fixed amount up or down in its vertical
geodesic, to obtain an isometry $\Phi \from Q_\gamma \to Q_{\gamma'}$. 

It follows that $D_{\gamma'}$ has an extra, nonisometric, affine symmetry,
and so $\gamma'$ is periodic.

This completes the proof of Lemma \ref{LemmaTeichPeriodic} and
therefore also of Proposition~\ref{proposition:periodicity:preserved}.


\section{The endgame: computing \protect{$\QI(\Gamma_H)$} and 
\protect{$\Comm(\Gamma_H)$}}
\label{section:endgame}

\subsection{An injection $\QI(X_H) \to \QSym(S^1)$}

Let $\QSym(S^1)$ denote the group of quasisymmetric homeomorphisms of the
circle $S^1$. It is well-known that the extension of a
quasi-isometry of $\hyp^2$ to $S^1$ induces an isomorphism $\QI(\hyp^2)
\approx \QSym(S^1)$: boundary extension of quasi-isometries defines an
injection from the quasi-isometry group of any word hyperbolic group to
the homeomorphism group of its boundary; and quasi-symmetric
homeomorphisms of $S^1$ are exactly the extensions of quasi-isometries of
$\hyp^2$. The key to all of this is the theorem of Ahlfors and Beurling
\cite{AhlforsBeurling} that quasisymmetric homeomorphisms of $S^1$ are
exactly the extensions of quasiconformal homeomorphisms of the unit disc.

Fix once and for all a base fiber $D_0=\pi^{-1}(\tau_0)$ of $X_H$,
$\tau_0\in T_H$. Identify $D_0$ with $\hyp^2$, so $\bdy D_0$ is
identified with $S^1$. Then there is a map $\Psi\from \QI(X_H)\to
\QSym(S^1)$, defined as follows. Given a quasi-isometry $f \from X_H
\to X_H$, the image $f(D_0)$ is Hausdorff close to some fiber $D'$.  
Consider the composition 
$$D_0 \to f(D_0) \to D' \to D_0
$$
where the first map is $f$ and
the other maps are closest point projections. Then the composition $D_0
\to D_0$ is clearly a quasi-isometry, and so induces an 
element $\QSym(S^1)$.  Hence we have a map $\Psi\from \QI(X_H)\to
\QSym(S^1)$, which is easily seen to be 
well-defined, and in fact a homomorphism.

\begin{proposition}
\label{PropQIinQSym}
The homomorphism $\Psi\from \QI(X_H)\to \QSym(S^1)$ is injective.
\end{proposition}

\begin{proof}
By Propositions \ref{proposition:hyperplanes:respected} and
\ref{proposition:horizontal:respected}, and the fact that the projection
$X_H \to T_H$ induces an isometry between the space of horizontal leaves
of $X_H$ with the Hausdorff metric and the quotient tree $T_H$, it follows
that any quasi-isometry $f \from X_H \to X_H$ induces a quasi-isometry of
$T_H$ and so $f$ induces a homeomorphism $\Xi(f) \from \Lambda \to
\Lambda$ of the Cantor set $\Lambda = \bdy T_H \subset\PMF$.

Suppose that $\Psi(f)\in \QSym(S^1)$ is the identity map on $S^1$.  

We claim that the induced map $\Xi(f)\from\Lambda\to \Lambda$ is the
identity. It follows that the induced quasi-isometry $T_H \to T_H$ is a
bounded distance from the identity, and so $f$ takes each horizontal leaf
of $X_H$ a bounded distance from itself; and since the induced boundary map
of that leaf is the identity, it follows that $f$ takes each point a
bounded distance from itself, proving the proposition.

For proving the claim (and for later purposes) we review the well-known
embedding of $\MF$ into the space of $\pi_1(\Sigma)$-invariant measures on
the ``\Mobius\ band beyond infinity'' of the hyperbolic plane. That is,
consider the \emph{double set} of the circle, $DS^1 = \{\{x,y\} \subset S^1
\suchthat x \ne y\}$; with respect to the Klein model $D_0 = \hyp^2 \subset
\R\P^2$, the usual duality gives a bijection between $DS^1$ and the \Mobius\
band beyond infinity $\R\P^2 - \overline \hyp^2$. 

Let $\M(DS^1)$ be the
space of Borel measures on $DS^1$ with the weak$^*$ topology, and let
$\PM(DS^1)$ be the space of projective classes of elements of $\M(DS^1)$.
The space $\MF$ embeds into the space of $\pi_1(\Sigma)$-invariant elements
of $\M(DS^1)$, by lifting a measured foliation on $\Sigma$ to a
$\pi_1(\Sigma)$-invariant measured foliation on $D_0$, and then identifying
each leaf of the lifted foliation with the correspond pair of endpoints in
$DS^1$. The space $\PMF$ therefore embeds as $\pi_1(\Sigma)$-invariant
elements of $\PM(DS^1)$. Given $b \in \PMF$ let $\mu_b \in \PM(DS^1)$ be
the corresponding projective class of measures, let $\Supp(\mu_b)
\subset DS^1$ be the support, and let $E_b=\abs{\Supp(\mu_b)}$ be the union
of all the pairs in $\Supp(\mu_b)$. Note that $E_b$ may also be described as
the \emph{endpoint set} of $b$, the set of endpoints in $S^1$ of the leaves
of the measured foliation $\wt\Fol_b$ on $D_0$ obtained by lifting any
measured foliation $\Fol_b$ on $\Sigma$ that represents $b$. 

By Proposition \ref{proposition:foliations:respected}, for all $b \in
\Lambda$ and all quasi-isometries $\phi$ of $X_H$ we have 
$$\Psi(\phi)(E_b)=E_{\Xi(\phi)(b)}
$$

From our assumption that $\Psi(f)$ is the identity on $S^1$ it follows that
$\Psi(f)(E_b)=E_b$, for all $b\in\Lambda$. Note however that if $b \ne b'
\in \Lambda$ then $E_b\intersect E_{b'} =\emptyset$, because the projective
classes of $\Fol_b,\Fol_{b'}$ in $\PMF$ are connected by a \Teichmuller\
geodesic in $\Hull\Lambda$, and so the measured foliations
$\Fol_b,\Fol_{b'}$ can be chosen to be transverse in
$\Sigma$; the lifted foliations $\wt\Fol_b, \wt\Fol_{b'}$ are
therefore transverse in $D_0$ and so their endpoint sets $E_b, E_{b'}$ in
$S^1$ are disjoint. It follows that $\Xi(f)(b)=b$ for all $b\in\Lambda$.
\end{proof}

\subsection{The orbifold $\O_H$ associated to a Schottky group $H$}
\label{SectionSubcover}

Our goal in this subsection and the next is to compute the quasi-isometry
group $\QI(\Gamma_H)$. Using the injection $\Psi\from \QI(\Gamma_H) \to
\QSym(S^1)$ provided by Proposition \ref{PropQIinQSym}, our computation
will consist of an explicit description of the subgroup
$\Psi(\QI(\Gamma_H)) \subgroup\QSym(S^1)$; without further mention we shall
identify $\QI(\Gamma_H)$ with this subgroup.

The first step in the computation of $\QI(\Gamma_H)$ is to find a
natural orbifold subcover $\Sigma \to \O_H$ associated to a Schottky
subgroup $H\subgroup \MCG(\Sigma)$. The orbifold $\O_H$ is an important
invariant of $H$; it is the smallest subcover of $\Sigma$ to which $H$
descends as a subgroup of $\MCG(\O_H)$.

Recall we have fixed a base fiber $D_0$ of $X_H$ over a base point
$\tau_0 \in T_H$, and we identify $D_0$ with the universal cover
$\wt\Sigma$.

To each periodic hyperplane $P_w$ we associate a finite-index supergroup of
$\pi_1(\Sigma)$ in $\QSym(S^1)$, as follows. In the associated singular
\solv-manifold $Q_{\gamma_w}$ pick any horizontal leaf $D_w$, which has an
induced $xy$ structure with $xy$-affine automorphism group
denoted $\Aff(D_w)$. Note that $\Aff(D_w)$ acts by quasi-isometries of the
$xy$-metric on $D_w$. There is a canonical quasi-isometry from $D_w$ to
$D_0$: restrict the canonical horizontal respecting quasi-isometry
$Q_{\gamma_w} \to P_w$ to the horizontal leaf $D_w$, giving a map $D_w \to
X_H$, and then take a closest point map to $D_0$. In this way we obtain an
inclusion
$$\Aff(D_w) \subset \QI(D_0) = \QSym(S^1)
$$
Note that the image of this inclusion is independent of the choice of
horizontal leaf $D_w$ in $Q_{\gamma_w}$, because if we chose another
horizontal leaf $D'_w$ then the vertical flow on $Q_{\gamma_w}$ induces an
$xy$-affine homeomorphism $D'_w \to D_w$ in the correct quasi-isometry
class.

Recall that we have a short exact sequence 
$$1 \to \Isom_+(D_w) \to \Aff(D_w) \to C_{\gamma_w} \to 1
$$ 
with $C_{\gamma_w} \equiv \Z$ or $D_\infinity$. Under the injection
$\Aff(D_w) \inject \QSym(S^1)$, the group $\Isom_+(D_w)$ is a finite index
supergroup of $\pi_1(\Sigma)$. We have a quotient orbifold $\O_w = D_w /
\Isom_+(D_w)$ with fundamental group $\pi_1(\O_w) \approx \Isom_+(D_w)$,
and associated to the inclusion $\pi_1(\O_w) \subset \pi_1(\Sigma)$ there
is a finite orbifold covering map $\Sigma \to \O_w$. 

In the group $\QSym(S^1)$, take the infinite intersection of the groups
$\Isom_+(D_w)$ over all periodic lines $w$ in the tree $T_H$, and note that
this group must be a finite index supergroup of $\pi_1(\Sigma)$, in fact it
is the fundamental group $\pi_1(\O_H)$ of an orbifold $\O_H$ which $\Sigma$
finitely covers, the smallest orbifold covered by $\Sigma$ which in turn
covers each orbifold $\O_w$.
$$\pi_1(\O_H) = \bigcap_{\substack{w \subset T_H \\ w
\text{ periodic}}}\Isom_+(D_w)
$$
Note that $\O_H$ is the smallest subcover of $\Sigma$ such that the
Schottky group $H \subset \MCG(\Sigma)$ descends via the covering map
$\Sigma \to \O_H$ to a free subgroup of $\MCG(\O_H)$. To be precise, since
$H$ is free we may choose a section $\sigma \from H \to \MCG(\Sigma,p)
\subset\QSym(S^1)$, and the image group $\sigma H$ acting by conjugation on
subgroups of $\QSym(S^1)$ permutes the collection of subgroups
$\Isom_+(D_w)$ and so $\sigma H \subset \Aut(\pi_1(\O_H)) = \MCG(\O_H,p)$.
Projecting to $\MCG(\O_H)$ we obtain the free subgroup $H'$, and a section
$\sigma' \from H' \to \MCG(\O_H,p) \subgroup \QSym(S^1)$ so that $\sigma H
= \sigma' H'$. 

\paragraph{Remark.} Orbifold mapping class groups obey many of the properties
of surface mapping class groups. The group $\MCG(\O_H)$ is defined as the
group of orbifold homeomorphisms of $\O_H$ modulo those which are isotopic
to the identity through orbifold homeomorphisms. Just as with mapping class
groups of surfaces, choosing a generic point $p \in \O_H$ we obtain an
isomorphism of short exact sequences
$$
\xymatrix{
1 \ar[r]  & \pi_1(\O_H) \ar[r] \ar@{=}[d]  & \MCG(\O_H,p) \ar[r]^\pi
\ar@2{~}[d]  &
\MCG(\O_H)
\ar[r] \ar@2{~}[d]  & 1 \\
1 \ar[r]  & \pi_1(\O_H) \ar[r]  & \Aut(\pi_1(\O_H)) \ar[r]^\pi  &
\Out(\pi_1(\O_H))  \ar[r]  & 1 }
$$
The proof that this isomorphism exists follows the same lines as the proof
for surfaces, using the fact that $\pi_1(\O_H)$ is centerless.

\bigskip

By choosing a lift of $p$ to the universal cover of $D_0 = \wt\O_H$ and
lifting each element of $\MCG(\O_H,p)$, we identify $\MCG(\O_H,p)$ with a
subgroup of $\QI(D_0) = \QSym(S^1)$. 

In fact $H' \subgroup \MCG(\O_H)$ is a Schottky subgroup. 
To see why, consider the limit set $\Lambda(H) \subset
\PMF(\Sigma)$. Moving into $S^1$, the corresponding subset $\{\mu_b
\suchthat b \in\Lambda(H)\}$ of $\PM(DS^1)$ is invariant under $\sigma H$.
Notice that for each $b \in\Lambda(H)$, the element $\mu_b \in \PM(S^1)$,
which is invariant under $\pi_1(\Sigma)$, is also invariant under the
larger group $\pi_1(\O_H)$. When $b$ is an endpoint of a periodic geodesic
$w$ in $T_H$ this is obvious, because $\mu_b$ is invariant under
$\pi_1(\O_w) \supset \pi_1(\O_H)$. But periodic endpoints are dense in
$\Lambda(H)$, and so each $\mu_b$, $b \in \Lambda(H)$ is invariant under
$\pi_1(\O_H)$. We therefore obtain a continuous embedding of
$\Lambda(H)$ in $\PMF(\O_H)$, whose image we denote $\Lambda(H')$; this
embedding has the property that
$$\{\mu_b \suchthat b \in \Lambda(H)\} = \{\mu_{b'} \suchthat b' \in
\Lambda(H')\} \quad\text{in}\quad \PM(DS^1)
$$
Moreover, since $\sigma' H' = \sigma H$, it follows that $\sigma' H'$
permutes the elements of the above set, and so $\Lambda(H')$ is in fact
invariant under the action of $H'$. For any pair $\xi\ne\eta
\in \Lambda(H)$ and corresponding pair $\xi',\eta' \in \Lambda(H')$, each
$xy$-structure on $\Sigma$ corresponding to a point along the
\Teichmuller\ geodesic $\geodesic{\xi}{\eta} \subset \T(\Sigma)$ lifts to
an $xy$-structures on $D_0$ which is invariant under $\pi_1(\O_H)$, and we
therefore obtain a  \Teichmuller\ geodesic $\geodesic{\xi'}{\eta'}$ in
$\T(\O_H)$; the $\pi_1(\O_H)$-invariance follows by an argument similar to
the one just above where we showed invariance for elements $\mu_b$, $b \in
\Lambda(H)$. The union of these geodesics over all pairs
$\xi' \ne \eta' \in
\Lambda(H')$ is denoted $\Hull\Lambda(H')$ as usual. From this construction
we see that the action of $H'$ on $\Lambda(H') \union \Hull\Lambda(H')$
agrees with the action of $H$ on $\Lambda(H) \union \Hull\Lambda(H)$. In
particular, all the requirements in Theorem \ref{theorem:schottkydef} for
$H'$ to be convex cocompact with limit set $\Lambda(H')$ are satisfied, and
so $H'$ is a Schottky subgroup of $\MCG(\O_H)$.

\paragraph{Remark.} The fact that the orbifold $\O_H$ supports a
pseudo-Anosov mapping class puts some restrictions on its topology: the
underlying surface of $\O_H$ must have empty boundary, because otherwise
the collection of peripheral curves would be invariant under any mapping
class, violating the existence of pseudo-Anosov mapping classes. It follows
that $\O_H$ is a closed surface with cone singularities.

\subsection{Computing the quasi-isometry group}
\label{SectionComputingQIGroup}

In this section we compute $\QI(\Gamma_H)$. Let $\C$ denote the relative
commensurator of $H'$ in $\MCG(\O_H)$:
$$\C = \Comm_{\MCG(\O_H)}(H')
$$
Form the extension group $\Gamma_\C$, a subgroup of
$\MCG(\O_H,p)$ as the following diagram shows:
$$
\xymatrix{
1 \ar[r]  & \pi_1(\O_H) \ar[r] \ar@{=}[d]  &
\Gamma_{\C} \ar[r]
\ar@{^{(}->}[d]  & \C \ar[r] \ar@{^{(}->}[d]  & 1 \\
1 \ar[r]  & \pi_1(\O_H) \ar[r]  & \MCG(\O_H,p) \ar[r]^\pi  &
\MCG(\O_H)  \ar[r]  & 1 
}
$$
The group $\Gamma_{\C}$ may therefore be regarded as a
subgroup of
$\QSym(S^1)$.

Here is our computation of $\QI(\Gamma_H)$:

\begin{theorem}
\label{TheoremQIGroup} 
In $\QSym(S^1)$ we have
$$\QI(\Gamma_H) = \Gamma_{\C}.
$$
Moreover, the subgroup $\pi_1(\O_H)$ consists of those classes of
quasi-isometries of $X_H$ which coarsely preserve each horizontal leaf of
$X_H$.
\end{theorem}

\begin{proof}
To justify this computation, first we show $\QI(\Gamma_H) \subset
\Gamma_{\C}$. Consider a quasi-isometry $f\from X_H\to
X_H$, regarded as an element of $\QSym(S^1)$; we must show that $f \in
\Gamma_{\C}$.

By Lemma \ref{PropABCCoarseRespect} the quasi-isometry $f$ coarsely
respects the horizontal foliation of $X_H$ and so $f$ lies over a
quasi-isometry of $T_H$, also denoted $f$. From Proposition
\ref{proposition:periodicity:preserved} it follows that $f$ induces a
permutation on the set of periodic hyperplanes of $X$: given a periodic
geodesic $w$ in $T_H$, if we let $w'$ be the geodesic in $T_H$ coarsely
equivalent to $f(w)$, then $f$ takes $P_w$ to $P_{w'}$ and
$Q_{\gamma_w}$ to $Q_{\gamma_{w'}}$. It follows that the conjugation action
of $f$ on subgroups of $\QSym(S^1)$ takes $\Isom(Q_{\gamma_w}) = \Aff(D_w)$
to $\Isom(Q_{\gamma_{w'}})=\Aff(D_{w'})$. Any isometry of $Q_{\gamma_w}$
which preserves each leaf of the horizontal foliation is conjugated by $f$
to a similar isometry of $Q_{\gamma_{w'}}$, and so conjugation by $f$ takes
$\Isom_h(Q_{\gamma_w}) = \Isom_+(D_w)$ to $\Isom_h(Q_{\gamma_{w'}}) =
\Isom_+(D_{w'})$. In other words, conjugation by $f$ preserves the
collection of subgroups 
$$\{\Isom_+(D_w) \suchthat  w \subset T_H \,\,\text{is periodic}\} 
$$
Intersecting this collection it follows that conjugation by $f$ in the
group $\QSym(S^1)$ preserves the subgroup $\pi_1(\O_H)$, acting as an
automorphism of that subgroup. In other words, in $\QSym(S^1)$ we have $f
\in \Aut(\pi_1(\O_H))=\MCG(\O_H,p)$.

Now we show that the image $\pi f \in \MCG(\O_H)$ lies in
$\C=\Comm_{\MCG(\O_H)}(H')$. We need the following fact:

\begin{theorem}[Commensurators of Schottky groups]
\label{theorem:normalizers}
Let $\O$ be a closed orbifold and $H'$ a Schottky subgroup of
$\MCG(\O)$.  Then the relative commensurator $\Comm_{\MCG(\O)}(H')$ of
$H'$ in $\MCG(\O)$ is equal to the subgroup of $\MCG(\O)$ stabilizing
the limit set $\Lambda(H')$, and $H'$ has finite index in
$\Comm_{\MCG(\O)}(H')$.
\end{theorem}

\begin{proof} Let $I$ be the subgroup of $\MCG(\O)$ stabilizing
$\Lambda(H')$, and so $H' \subgroup I$. Let $\Hull\Lambda(H')$ be the weak
convex hull of $\Lambda(H')$ (see Theorem \ref{theorem:schottkydef}).  
Clearly $I$ is also the subgroup of $\MCG(\O)$ stabilizing
$\Hull\Lambda(H')$. By Theorem \ref{theorem:schottkydef} the group $H'$ acts
cocompactly on $\Hull\Lambda(H')$, and so $I$ acts cocompactly on
$\Hull\Lambda(H')$. Clearly $I$ acts properly on $\Hull\Lambda(H')$. It
follows that $H'$ has finite index in $I$, which immediately implies $I
\subgroup \Comm_{\MCG(\O)}(H')$.

For the opposite inclusion, suppose $\Phi \in \MCG(\O) - I$, and so
$\Phi(\Lambda(H')) \ne \Lambda(H')$. The limit set of the Schottky subgroup
$\Phi H' \Phi^\inv$ is $\Phi(\Lambda(H'))$. Since $\Lambda(H')$ is closed,
$\Phi(\Lambda(H')) - \Lambda(H')$ is open in $\Phi(\Lambda(H'))$. Since
fixed points of pseudo-Anosov elements of $\Phi H' \Phi^\inv$ are dense in
$\Phi(\Lambda(H'))$, there exists a pseudo-Anosov element $\Psi \in \Phi H'
\Phi^\inv$ having a fixed point not in $\Lambda(H')$. Infinitely many
powers of $\Psi$ are therefore not in $H'$, and so $H' \intersect \Phi H'
\Phi^\inv$ has infinite index in $\Phi H' \Phi^\inv$. \end{proof}

We showed above that $H'$ is a Schottky subgroup of $\MCG(\O_H)$ and so 
by Theorem \ref{theorem:normalizers} it
remains to show that $\pi f \in \MCG(\O_H)$ acting on $\PMF(\O_H)$
leaves the set $\Lambda(H')$ invariant. For this purpose it suffices to
show that the action of $f \in \QSym(S^1)$ leaves the
set $\{\mu_{b'} \suchthat b' \in \Lambda(H')\}$ invariant. However, we know
that this set equals $\{\mu_b \suchthat b \in \Lambda(H)\}$, and we also
know that $f$ leaves this set invariant: $f$ permutes the elements
of this set which are endpoints of periodic geodesics in $T_H$, by
Proposition \ref{proposition:periodicity:preserved}, but the endpoints of
periodic geodesics are dense. 

This completes the proof of the inclusion $\QI(\Gamma_H) \subset
\Gamma_{\C}$. 

For the opposite inclusion, consider the following commutative diagram of
short exact sequences:
\begin{equation}
\label{EquationSESDiagram}
\xymatrix{
1 \ar[r] & \pi_1(\Sigma) \ar[r] \ar@{^{(}->}[d] & \Gamma_H \ar[r]
\ar[d] & H \ar[r] \ar@{^{(}->}[d] & 1 \\
 1 \ar[r]  & \pi_1(\O_H) \ar[r] & \Gamma_{\C} \ar[r]
 & \C \ar[r] & 1 
}
\end{equation}
The left vertical arrow is a finite index injection because $\Sigma \to
\O_H$ is a finite covering. The right vertical arrow is a finite index
injection by Theorem \ref{theorem:normalizers} using the isomorphism $H
\approx H'$. It follows that the middle vertical arrow is a finite index
injection. But this shows that $\Gamma_{\C}$ is a finite index
supergroup of $\Gamma_H$ in $\QSym(S^1)$, and so $\Gamma_{\C} \subset
\QI(\Gamma_H)$.

This completes the computation of $\QI(\Gamma_H)$. 

\bigskip

To complete the proof of Theorem~\ref{TheoremQIGroup},
it remains to prove that an element of $\QI(\Gamma_H) = \Gamma_{\C}$ is
in the normal subgroup $\pi_1(\O_H)$ if and only if it coarsely preserves
each horizontal leaf of
$X_H$. Evidently $\pi_1(\O_H)$ coarsely preserves each horizontal leaf.
Suppose conversely that the quasi-isometry $f \from X_H \to X_H$ coarsely
preserves each horizontal leaf; let $[f] \in \QI(\Gamma_H)$ be the
coarse equivalence class of
$f$. It follows that $f$ coarsely respects each hyperplane of
$X_H$, in particular each periodic hyperplane $P_w$, and so we have $[f]
\in \Aff(D_w)$ for each periodic line $w$ in $T_H$. But we can say more, as
follows. The inclusion map $\Gamma_H \to \QI(\Gamma_H)$, being injective
and with finite index image, is a quasi-isometry, and we therefore obtain a
quasi-isometry $\QI(\Gamma_H) \to X_H$. Using this quasi-isometry we may
conjugate the left action of $\QI(\Gamma_H)$ on itself to obtain a
quasi-action of $\QI(\Gamma_H)$ on $X_H$. In particular, we obtain a
sequence of uniform quasi-isometries $f_n \from X_H \to X_H$, with $[f_n] =
[f]^n$. Each of the $f_n$ coarsely preserves each
horizontal leaf of $X_H$, and note that the coarseness constant is uniform
\emph{independent of $n$ and of the leaf}, by application of Proposition
\ref{proposition:hyperplanes:respected} using uniformity of the
quasi-isometry constants of $f_n$. It follows that each $f_n$ coarsely
preserves each periodic hyperplane $P_w$, and $f_n$ coarsely preserves each
horizontal leaf of $P_w$, again with uniform coarseness constants
independent of $n$. This implies that $[f] \in \Isom_+(D_w)$. We
therefore have
$$[f] \in \bigcap_w \Isom_+(D_w) = \pi_1(\O_H)
$$
\end{proof}

\subsection{Computing the commensurator group}
\label{section:comm}

In this section we use our knowledge of quasi-isometries of $\Gamma_H$ to
compute the commensurator of $\Gamma_H$ and complete the proof of
Theorem~\ref{TheoremCommQI}. This is the first instance we
know of where this technique of computing a commensurator group is used. 
We prove, in Theorem~\ref{TheoremCommQIIsomorphism}, that when the left
multiplication homomorphism $\Gamma \to \QI(\Gamma)$ is an injection with
finite index image, the natural homomorphism $i \from \Comm(\Gamma) \to
\QI(\Gamma)$ is an isomorphism. Combining with
Theorem~\ref{TheoremQIGroup} we obtain the desired computation of
$\Comm(\Gamma_H)$. We do not have a proof of the computation
of $\Comm(\Gamma_H)$ without going through all of the work required to
understand quasi-isometries; a purely algebraic computation of
$\Comm(\Gamma_H)$ would be interesting. 

For any finitely-generated group $\Gamma$, a commensuration $\phi\from
G_1\to G_2$ extends to a quasi-isometry of $\Gamma$ by precomposing
$\phi$ with a closest point projection $\Gamma\to G_1$ and postcomposing
with inclusion $G_2 \to \Gamma$. The coarse equivalence class of this
extension is well-defined, giving a natural homomorphism
$i\from\Comm(\Gamma)\to\QI(\Gamma)$. Note that we have a commutative
triangle
$$
\xymatrix{
\Gamma \ar[rrr]^{C \from g \mapsto \{C_g\}} \ar@/^2pc/[rrrrr]^{L \from g
\mapsto [L_g]} & & &
\Comm(\Gamma)
\ar[rr]^i & &
\QI(\Gamma) }
$$
where $L_g$ is the left multiplication map $x \to gx$ with coarse
equivalence class $[L_g]\in\QI(\Gamma)$,
and
$C_g$ is the conjugation automorphism $x \to gxg^\inv$ with equivalence
class $\{C_g\}\in\Comm(\Gamma)$; commutativity follows because right
multiplication by $g^\inv$ moves each point in $\Gamma$ a bounded amount.

\begin{theorem}
\label{TheoremCommQIIsomorphism}
Given a finitely generated group $\Gamma$, if $L\from\Gamma\to
\QI(\Gamma)$ is an injection with finite index image, then $i \from
\Comm(\Gamma)\to\QI(\Gamma)$ is an isomorphism.
\end{theorem}

This theorem, proved in collaboration Kevin Whyte, is broken into three
steps: 
\begin{enumerate}
\item injectivity of $i\from \Comm(\Gamma)\to\QI(\Gamma)$; 
\item construction of an injective map
$\Psi\from\QI(\Gamma)\to\Comm(\Gamma)$;
\item the proof that $\Psi\composed i$ is the identity on
$\Comm(\Gamma)$. 
\end{enumerate}

\paragraph{Step 1: Injectivity of $i$.} 

\begin{proposition}[K.\ Whyte]
\label{prop:CommQIInject}
If $\Gamma$ is any finitely generated group then the natural map
$i \from \Comm(\Gamma)\to \QI(\Gamma)$ is injective.
\end{proposition}

\begin{proof}
Suppose $\phi\from G_1\to G_2$ is a commensuration of
$\Gamma$ such that $i(\phi)$ equals the identity in $\QI(\Gamma)$. 
Then there is a bounded function $\delta\from\Gamma\to \Gamma$
so that for all $g\in G_1$ we have
$$\phi(g)=g \, \delta(g)
$$
Since $\delta$ is bounded, the cardinality $M = \#\image(\delta)$ is
finite.
 
Plugging the above equation into $\phi(gh)=\phi(g)\phi(h)$ 
gives 
$$h^{-1}\delta(g)h=\delta(gh)\delta(h)^\inv
$$
Note that this is true for all $h \in G_1$, and the right hand side takes
on at most $M^2$ values. This implies that the centralizer of $\delta(g)$
in $G_1$ has index at most $M^2$. Intersecting all subgroups of $G_1$ of
index $\le M^2$ gives a finite index subgroup $H \subgroup G_1$ which
commutes with each $\delta(g)$. From the above equation it follows that
$\delta$ is a homomorphism on $H$. Since $\image(\delta)$ is finite it
follows that $\delta$ has finite index kernel $\kernel(\delta)$ in
$H$, and so $\kernel(\delta)$ has finite index in $G_1$ and in
$\Gamma$. In other words,
$\phi(g)=g$ on the finite index subgroup $\kernel(\delta)$ of $\Gamma$,
and so $\phi$ represents the identity element of $\Comm(\Gamma)$.
\end{proof}

For any finitely generated group $\Gamma$ the kernel of $C \from \Gamma
\to \Comm(\Gamma)$ is the \emph{virtual center} $\VZ(\Gamma)$ consisting
of all $g \in \Gamma$ whose centralizer has finite index in $\Gamma$.
Together with Proposition \ref{prop:CommQIInject} it follows that
injectivity of $L \from \Gamma\to\QI(\Gamma)$ is equivalent to the
triviality of $\VZ(\Gamma)$. 

\paragraph{Step 2: An injection $\Psi \from \QI(\Gamma) \to
\Comm(\Gamma)$.} 

Identifying $\Gamma$ with its image $L(\Gamma) \subgroup \QI(\Gamma)$, a
finite index subgroup, it follows that any automorphism of $\QI(\Gamma)$
restricts to a commensuration of $\Gamma$. In particular, given $F \in
\QI(\Gamma)$, the inner automorphism $G \mapsto FGF^\inv$ of
$\QI(\Gamma)$ restricts to a commensuration $\Psi_F$ of $\Gamma$, giving a
well-defined homomorphism $\Psi\from \QI(\Gamma) \to \Comm(\Gamma)$.
Keeping in mind the definition of $L \from \Gamma \to \QI(\Gamma)$, this
means that for each $F\in\QI(\Gamma)$ and each $x \in \Gamma$ we have the
following equation in $\QI(\Gamma)$:
$$
F \cdot [L_x] \cdot F^\inv = [L_{\Psi_F(x)}]
$$

To prove that $\Psi$ is an injection, given $F=[f] \in \QI(\Gamma)$
suppose that $\Psi_F$ is the identity map when restricted to a finite
index subgroup $G_1$ of $\Gamma$. Then for all $x \in G_1$ we have
\begin{align*}
F \cdot [L_x] \cdot F^\inv &= [L_x] \\
F \cdot [L_x]  &= [L_x] \cdot F \\
[f \composed L_x] &= [L_x \composed f]
\end{align*}
which means that there exists a bounded function $\delta \from \Gamma \to
\Gamma$ such that
$$f(x y)  = x f(y) \delta(y), \quad x, y \in G_1
$$
Plugging in $y=1$ we get
$$f(x) = x f(1) \delta(1), \quad x \in G_1
$$
This shows that the function $f$ is a bounded distance from the identity
map on $G_1$, and so $F$ is the identity element in $\QI(\Gamma)$.

\paragraph{Step 3: $\Psi \composed i \from \Comm(\Gamma)
\to \Comm(\Gamma)$ is the identity.} Given a commensuration $\phi$ of
$\Gamma$, we obtain a commensuration $\phi'=\Psi_{i(\phi)}$ satisfying the
formulas
\begin{align*}
[\phi \composed L_x \composed \phi^\inv] &= [L_{\phi'(x)}] \\
[\phi \composed L_x \composed \phi^\inv \composed L_{\phi'(x^\inv)}] &=
[\Id]
\end{align*}
for all $x$ in a certain finite index subgroup of $\Gamma$. We must prove
that $\phi$ and $\phi'$ agree on a further finite index subgroup of
$\Gamma$.

From the above equation it follows that there is a bounded function
$\delta \from \Gamma \to \Gamma$ such that for all $h$ in a certain
finite index subgroup $H$ of $\Gamma$ we have:
\begin{align*}
\phi(x \cdot \phi^\inv(\phi'(x^\inv) \cdot h) &= h \cdot \delta(h) \\
\phi(x) \cdot \phi'(x^\inv) \cdot h &= h \cdot \delta(h) \\
h^\inv  \cdot \bigl( \phi(x) \cdot \phi'(x^\inv) \bigr) \cdot h &=
\delta(h) 
\end{align*}
This shows that $\phi(x) \cdot \phi'(x^\inv)$ has only finitely many
conjugates by elements of $H$, and so the centralizer of $\phi(x) \cdot
\phi'(x^\inv)$ in $H$ has finite index in $H$. It follows that the
centralizer of $\phi(x) \cdot \phi'(x^\inv)$ in $\Gamma$ has finite index
in $\Gamma$, which by definition of $\VZ(\Gamma)$ gives $\phi(x) \cdot
\phi'(x^\inv)\in\VZ(\Gamma)$. But since $L \from \Gamma \to \QI(\Gamma)$
is injective, the virtual center $\VZ(\Gamma)$ is trivial, showing that
$\phi(x) = \phi'(x)$ for all $x$ in a finite index subgroup of $\Gamma$.

\medskip

Combining steps 1--3 it follows that $\Psi \from \QI(\Gamma) \to
\Comm(\Gamma)$ is surjective, and so $\Psi$ is an isomorphism with
inverse $i \from \Comm(\Gamma) \to \QI(\Gamma)$, finishing the proof of
Theorem~\ref{TheoremCommQIIsomorphism}.

\begin{Remark}
The proof of Theorem~\ref{TheoremCommQIIsomorphism} yields a more general
conclusion: if $\Gamma$ is a finitely generated group whose virtual
center is trivial, then $\Comm(\Gamma)$ is isomorphic to the relative
commensurator of $L(\Gamma)$ in $\QI(\Gamma)$. The condition that the
virtual center be trivial cannot be dropped: for example, a finite group
has trivial quasi-isometry group but rarely is its abstract commensurator
group trivial.
\end{Remark}

\subsection{Proving Theorem \ref{theorem:classification} and Theorem
\ref{theorem:rigidity}}
\label{SectionClassAndRigidityProofs}

Theorem \ref{theorem:rigidity} follows immediately from 
Theorem \ref{theorem:finiteindex} by a standard technique (see,
e.g.\ \cite{Schwartz:RankOne}). The basic observation we need says that if
$\Gamma$ is a finitely generated group and if the homomorphism $\Gamma \to
\QI(\Gamma)$ has finite cokernel and kernel, then for any finitely
generated group $H$ and any quasi-isometry $\phi \from H \to \Gamma$ the
induced homomorphism $\phi_* \from H \to \QI(\Gamma)$ has finite index
kernel and cokernel. As a consequence, the groups $H$ and $\Gamma$ are
weakly commensurable, because their images are commensurable in
$\QI(\Gamma)$. 

We now prove Theorem \ref{theorem:classification}. Let $H_i \subgroup
\MCG(\Sigma_{g_i})$, $i=1,2$, be Schottky groups of rank $\ge 2$ with $g_i
\ge 2$. We must prove the equivalence of the following four statements:
\begin{itemize}
\item[(\ref{ItemQI})]
$\Gamma_{H_1}$ and $\Gamma_{H_2}$ are quasi-isometric.
\item[(\ref{ItemAbstractComm})]
$\Gamma_{H_1}$ and $\Gamma_{H_2}$ are abstractly commensurable. 
\item[(\ref{ItemConcreteComm})]
There is an isomorphism $\O_{H_1} \approx \O_{H_2}$ such that in the group
$\MCG(\O_{H_1}) = \MCG(\O_{H_2})$ the Schottky subgroups $H_1$ and $H_2$ are
commensurable, meaning that $H_1 \intersect H_2$ has finite index in each
of $H_1$ and $H_2$.
\item[(\ref{ItemSameLimitSet})]
There is an isomorphism $\O_1 \approx \O_2$ such that in the group
$\MCG(\O_1) = \MCG(\O_2)$ the Schottky groups $H_1$ and $H_2$ have the same
limit set in the Thurston boundary of the \Teichmuller\ space $\Teich(\O_1)
= \Teich(\O_2)$.
\end{itemize}
We also add in a fifth equivalent statement:
\begin{itemize}
\item[(5)] The groups $\QI(\Gamma_{H_1})$ and $\QI(\Gamma_{H_2})$ are
isomorphic.
\end{itemize}
The equivalence of statements (\ref{ItemQI}) and (\ref{ItemAbstractComm})
and (5) follows immediately from Theorem \ref{theorem:finiteindex} using the
fact that a quasi-isometry between two groups induces an isomorphism between
their quasi-isometry groups. 

The fact that (\ref{ItemConcreteComm}) implies (\ref{ItemAbstractComm}) is
an immediate consequence of the commutative diagram
\ref{EquationSESDiagram} applied to
$\Gamma_{H_1}$ and to $\Gamma_{H_2}$. 

To prove that (\ref{ItemQI}) implies (\ref{ItemConcreteComm}),
suppose that there is a quasi-isometry $\Gamma_{H_1} \to \Gamma_{H_2}$,
which induces an isomorphism of quasi-isometry groups
$\QI(\Gamma_{H_1}) \approx \QI(\Gamma_{H_2})$.  Consider, for each $i=1,2$,
the short exact sequence
$$1 \to \pi_1(\O_i) \to \QI(\Gamma_{H_i}) \to \C_i \to 1
$$
where as usual $\O_i = \O_{H_i}$ is the smallest orbifold subcover of
$\Sigma_{g_i}$ to which $H_i$ descends, the subgroup $\C_i \subgroup
\MCG(\O_i)$ is the relative commensurator of $H_i$ in $\MCG(\O_i)$, and
$H_i$ has finite index in $\C_i$.

We claim that the isomorphism between $\QI(\Gamma_{H_1})$ and
$\QI(\Gamma_{H_2})$ must take $\pi_1(\O_1)$ to $\pi_1(\O_2)$. This provides
an isomorphism $\O_1 \approx \O_2$ such that the induced isomorphism
$\MCG(\O_1) \approx \MCG(\O_2)$ takes $\C_1$ to $\C_2$, and statement
(\ref{ItemConcreteComm}) immediately follows. 

To prove the claim, consider the model space $X_{H_i}$ with its horizontal
foliation. By Theorem \ref{TheoremQIGroup} the subgroup $\pi_1(\O_i)$ of
$\QI(\Gamma_{H_i})\approx \QI(X_{H_i})$ consists of quasi-isometries which
coarsely preserve each leaf of the horizontal foliation. The quasi-isometry
$X_{H_1} \to X_{H_2}$ coarsely respects horizontal foliations by Lemma 
\ref{PropABCCoarseRespect}, and so the coarse leaf preserving
elements of $\QI(X_{H_1})$ are taken bijectively by the isomorphism
$\QI(X_{H_1}) \leftrightarrow \QI(X_{H_2})$ to the coarse leaf respecting
quasi-isometries of $X_{H_2}$. In other words, $\pi_1(\O_1)$ is taken to
$\pi_1(\O_2)$.

Finally, the equivalence of (\ref{ItemConcreteComm}) and
(\ref{ItemSameLimitSet}) follows immediately from
Theorem~\ref{theorem:normalizers}, completing the proof of Theorem
\ref{theorem:classification}.
\endproof

\section{Closing remarks}

\subsection{Surfaces versus orbifolds}
\label{SectionSurfaceVsOrbifold}

As remarked in the introduction, the
universe of groups quasi-isometric to (orbifold)-by-(virtual Schottky)
groups is exactly the same as the universe of groups quasi-isometric to
(surface)-by-(Schottky) groups:

\begin{proposition} 
\label{PropSurfaceVsOrbifold}
Given a closed orbifold $\O$ and a virtual Schottky subgroup $N \subgroup
\MCG(\O)$, consider the extension $\Gamma_N$ defined by 
$$1 \to \pi_1(\O) \to \Gamma_N\to N \to 1 
$$ 
There exists a closed oriented
surface $\Sigma$ of genus $g \ge 2$, and a Schottky subgroup $H \subgroup
\MCG(\Sigma)$, such that the group $\Gamma_H = \pi_1(\Sigma) \semidirect H$
has finite index in $\Gamma_N$, and so the two groups are quasi-isometric.
\end{proposition}

\begin{proof}
There is a Schottky subgroup $H' \subgroup N$ of finite index; it follows
that $\Gamma_{H'}$ has finite index in $\Gamma_N$. Choose a splitting $H'
\to\MCG(\O,p)$, consider the action of $H'$ on $\pi_1(\O)$ by
automorphisms. Choose a finite surface cover $\Sigma'$ with corresponding
subgroup $\pi_1(\Sigma')
\subgroup \pi_1(\O)$, and consider the orbit of $\pi_1(\Sigma')$ under the
action of $H'$. This orbit consists of a finite collection of finite index
subgroups of $\pi_1(\O)$, whose intersection is a finite index subgroup
corresponding to a surface group $\Sigma$ which is a cover of $\Sigma'$.
The group $H' \subgroup \MCG(\O,p)$ lifts to a subgroup $H \subgroup
\MCG(\Sigma,p)$, which projects to a Schottky subgroup of
$\MCG(\Sigma)$. The group $\Gamma_{H}=\pi_1(\Sigma) \semidirect H$ is
therefore a (surface)-by-(Schottky) group with finite index in
$\Gamma_{H'}$, and so also in the original group
$\Gamma_N$.
\end{proof}

\subsection{Fibered hyperbolic 3-manifold groups}
\label{SectionHorResClasses}

As mentioned in the footnote on page~\pageref{Footnote}, the method of
proof of Theorem \ref{theorem:classification} shows that if two word
hyperbolic surface-by-free groups are quasi-isometric then they are
horizontal  respecting quasi-isometric, indeed they are horizontal
respecting commensurable, as long as the free group has rank $\ge 2$.
When the free group is infinite cyclic this fails: all
hyperbolic 3-manifolds fibering over the circle have quasi-isometric
fundamental groups; but there exist fibered hyperbolic 3-manifolds which
are not abstract commensurable---take, for example, an arithmetic example
and a non-arithmetic example.

Nevertheless, we do obtain a classification of fundamental groups of
fibered hyperbolic 3-manifold groups up to horizontal respecting
quasi-isometry. Namely, let $\Sigma$ be a closed surface of genus $g \ge
2$, let $\psi \in \MCG(\Sigma)$ be a pseudo-Anosov mapping class
generating an infinite cyclic subgroup $\generatedby{\psi}$ of
$\MCG(\Sigma)$, and let $\Gamma_{\generatedby{\psi}}$ be the associated
extension group of $\pi_1(\Sigma)$ by $\generatedby{\psi}$. Let
$\O_\generatedby{\psi}$ be the smallest subcover of $\Sigma$ to which
$\generatedby{\psi}$ descends, let $\C_\generatedby{\psi}$ be the relative
commensurator of $\generatedby{\psi}$ in $\MCG(\O_\generatedby{\psi})$,
and let $\Gamma_{\C_\generatedby{\psi}}$ be the extension of
$\pi_1(\O_\generatedby{\psi})$ by $\C_\generatedby{\psi}$. Note that
$\C_{\generatedby{\psi}}$ contains $\generatedby{\psi}$ with finite index
by \cite{BLMcC}, and so $\Gamma_{C_{\generatedby{\psi}}}$ contains
$\Gamma_{\generatedby{\psi}}$ with finite index. The proofs of
\S\ref{section:endgame} now apply directly to obtain the following
classification theorem, which parallels
Theorem~\ref{theorem:classification}:

\begin{theorem} Given closed surfaces $\Sigma_i$, $i=1,2$, and
pseudo-Anosov mapping classes $\psi_i$ of $\MCG(\Sigma_i)$, the
following are equivalent, where $\O_i = \O_{\generatedby{\psi_i}}$, etc.:
\begin{enumerate}
\item There exists
a horizontal respecting quasi-isometry $\Gamma_{\generatedby{\psi_1}} \to
\Gamma_{\generatedby{\psi_2}}$.
\item The groups $\Gamma_{\generatedby{\psi_1}},
\Gamma_{\generatedby{\psi_2}}$ are horizontal respecting commensurable,
meaning that there is an isomorphism from a finite index subgroup $H_1
\subgroup \Gamma_{\generatedby{\psi_1}}$ to a finite index subgroup $H_2
\subgroup
\Gamma_{\generatedby{\psi_2}}$, taking $H_1
\intersect \pi_1(\O_1)$ to $H_2 \intersect \pi_1(\O_2)$. 
\item
There is an isomorphism $\O_1 \approx \O_2$ such that in the group
$\MCG(\O_1) = \MCG(\O_2)$ the mapping classes $\psi_1$ and
$\psi_2$ have equal powers.  
\item
There is an isomorphism $\O_1 \approx \O_2$ such that in the group
$\MCG(\O_1) = \MCG(\O_2)$ the mapping classes $\psi_1$ and $\psi_2$ have
the same limit set in the Thurston boundary of the \Teichmuller\ space
$\Teich(\O_1) = \Teich(\O_2)$, i.e.\ they have the same stable/unstable
measured foliation pairs.
\item There is an isomorphism $\Gamma_{\C_1} \approx
\Gamma_{\C_2}$ taking $\pi_1(\O_1)$ to $\pi_1(\O_2)$.
\end{enumerate}
\qed\end{theorem} 

It is easy to use this theorem to obtain infinitely many distinct
horizontal respecting quasi-isometry classes of groups $\Gamma_H$, for
cyclic pseudo-Anosov groups $H$. This should be contrasted with the fact
that \emph{all} of the groups $\Gamma_H$ are quasi-isometric to each other
and to $\hyp^3$, by Thurston's hyperbolization theorem for fibered
\nb{3}manifolds \cite{Otal:fibered}.


\newcommand{\etalchar}[1]{$^{#1}$}
\providecommand{\bysame}{\leavevmode\hbox to3em{\hrulefill}\thinspace}

\noindent
Benson Farb:\\ Department of Mathematics, Univeristy of Chicago\\ 5734
University Ave.\\ Chicago, Il 60637\\ E-mail: farb@math.uchicago.edu
\medskip

\noindent
Lee Mosher:\\
Department of Mathematics, Rutgers University, Newark\\
Newark, NJ 07102\\
E-mail: mosher@andromeda.rutgers.edu

\end{document}